\documentclass[12pt]{article}

\usepackage[english]{babel}
\usepackage{indentfirst}
\usepackage{latexsym}
\usepackage[colorlinks=true, bookmarksnumbered=true, bookmarksopen=true,
bookmarksopenlevel=3, pdfstartview=FitH, linkcolor=cyan, pdfmenubar=true,
pdftoolbar=true, bookmarks=true,citecolor=cyan, urlcolor=magenta,
filecolor=cyan,menucolor=black,plainpages=false,pdfpagelabels]{hyperref}
\usepackage[lmargin=2cm,tmargin=2cm,bmargin=2cm,rmargin=2cm]{geometry}
\usepackage{amsbsy, amsmath,amsfonts,amssymb,amsthm}
\usepackage{times}
\usepackage{graphicx}
\usepackage{amsmath}
\usepackage{amsfonts}
\usepackage{amssymb}
\usepackage[english]{babel}
\usepackage[latin1]{inputenc}
\usepackage[T1]{fontenc}
\usepackage{amsmath,amssymb,graphics,amsthm}
\usepackage[usenames]{color}

\numberwithin{equation}{section}
\newtheorem{prop}{Proposition}[section]
\newtheorem{defi}[prop]{Definition}
\newtheorem{teo}[prop]{Theorem}

\newtheorem{obs}[prop]{Remark}
\newtheorem{lema}[prop]{Lemma}

\newcommand\supp{\mathop{\rm supp}}
\newcommand\reallywidehat[1]{\savestack{\tmpbox}{\stretchto{  \scaleto{    \scalerel*[\widthof{\ensuremath{#1}}]{\kern-.6pt\bigwedge\kern-.6pt}    {\rule[-\textheight/2]{1ex}{\textheight}}  }{\textheight}}{0.5ex}}\stackon[1pt]{#1}{\tmpbox}}

\begin{document}

\title{Global well-posedness for the Navier-Stokes system 
\\
in new critical mixed-norm Besov spaces}
\author{
{ Leithold L. Aurazo-Alvarez$^{1}$}{\thanks{
{E-mail address: aurazoall@gmail.com}
\newline
{L.L. Aurazo-Alvarez was supported by FAPERJ (Grants 200.140/200.141/2024) and Federal University of Rio de Janeiro, RJ, Brazil.
 }}} , 
\, {Wladimir Neves $^{2}$}{\thanks{{Corresponding author. }
\newline
{E-mail address: wladimir@im.ufrj.br}\newline
{Wladimir Neves has received research grants from FAPERJ (Cientista do Nosso Estado) through the grant E-26/204.171/2024, 
and also from CNPq through the grant  313005/2023-0, 406460/2023-0.
 }}} 
 \\[5pt]
{\small $^{1,2}$ Instituto de Matem\'atica, Universidade Federal do Rio de Janeiro } 
\\[5pt]
{\small Cidade Universit\'aria, CEP 21945-970, Rio de Janeiro, Brazil.}
}
\date{}
\maketitle

\begin{abstract}	
In this work, we proved the existence of a unique global mild solution of the d-dimensional incompressible 
Navier-Stokes equations, for small initial data in Besov type spaces based on mixed-Lebesgue spaces; namely, 
mixed-norm Besov-Lebesgue spaces and also mixed-norm Fourier-Besov-Lebesgue spaces. The main tools are the Bernstein's type inequalities, 
Bony's paraproduct to estimate the bilinear term and a fixed point scheme in order to get the well-posedness. 
Our results complement and cover previous and recents result on (Fourier-)Besov spaces and, for instance, 
provide a new class of initial data possibly not included in $BMO^{-1}(\mathbb{R}^{3})$ but continuously 
included in $\dot{{\cal B}}^{-1}_{\infty,\infty}(\mathbb{R}^{3})$.  

 \medskip
 

{\small \medskip\noindent\textbf{Key:} 
Well-posedness; Navier-Stokes equations; Critical spaces; Besov spaces; Fourier-Besov spaces; Mixed-Lebesgue spaces; Large initial data.}
\end{abstract}

\renewcommand{\abstractname}{Abstract}

\section{Introduction}
Considering a new class of Besov type spaces, we study in this paper the $d$-dimensional incompressible Navier-Stokes equations, that is to say
\begin{equation}
\label{BQ}
\left\{
\begin{aligned}
&  \partial_{t} u(x,t) + (u(x,t) \cdot \nabla) u(x,t) 
+ \nabla P(x,t) = \nu \Delta u(x,t),
\\[5pt]
&  \mbox{div}\,u(x,t)= 0,\quad \mbox{for}\,\,x=(x_1, \dots, x_d) \in \mathbb{R}^{d}, \quad t> 0,
\\[5pt]
&  u(x,0)=u_{0}(x),\,\,\mbox{for}\,x\in\mathbb{R}^{d},
\end{aligned}
\right.  
\end{equation}
where $u=u(x,t)=(u^{1}(x,t), \dots, u^{d}(x,t))$ represents the velocity vector field, the scalar function $P=P(x,t)$ denotes the pressure, 
the parameter $\nu> 0$ represents the kinematic viscosity, and $u_{0}$ is a given divergence free initial velocity field. From now on we call this system as ``Navier-Stokes system'' or ``Navier-Stokes equations'', however it is clear that we do not consider the compressible case.  

\medskip
First, let us recall some previous works on Besov or Fourier-Besov spaces, mostly dealing with global 
results for the Navier-Stokes system. Related to Fourier-Besov spaces, Z. 
Lei and F. H. Lin \cite{Lein-Lin-2011} established the global well-posedness for the $3D$-Navier-Stokes equations for small initial data 
in the now called Lei-Lin class ${\cal X}^{-1}$ corresponding to the space $\dot{{\cal FB}}^{-1}_{1,1}(\mathbb{R}^{3})$. 
Also, in the work \cite{Konieczny-Yoneda-2011},  P. Konieczny and T. Yoneda studied the dispersive effect of the Coriolis 
force for the non-stationary Navier-Stokes equations as well as for the stationary Navier-Stokes system. In the first part of 
that paper, they obtain the global well-posedness for the tridimensional Navier-Stokes system, for instance, for small initial 
data in the space $\dot{{\cal FB}}^{2-\frac{3}{p}}_{p,\infty}(\mathbf{R}^{3})$, for $1<p\leq \infty$, in the space 
$\dot{{\cal FB}}^{2-\frac{3}{p}}_{p,p}(\mathbf{R}^{3})$, for $3<p< \infty$, or in the space $\dot{{\cal FB}}^{-1}_{1,1}(\mathbf{R}^{3})\cap 
\dot{{\cal FB}}^{0}_{1,1}(\mathbf{R}^{3})$. These results are uniform in relation to the Coriolis parameter, in the sense that, the smallness condition does not depends on the Coriolis parameter. At the second part of that work, they analize the stationary Navier-Stokes system, and make use of the one-norm approach in the nonstationary case, for instance, they show that, for an given external force  $F\in \dot{{\cal FB}}^{-\frac{3}{p}}_{p,p}(\mathbf{R}^{3})$, with $3<p<\infty$, it is possible to find a Coriolis parameter $\Omega_{0}$ such that, for all $\mid\Omega\mid\geq \Omega_{0}$ there is a unique solution $u\in \dot{{\cal FB}}^{2-\frac{3}{p}}_{p,p}(\mathbf{R}^{3})$ for the stationary system. 

\medskip
To follow, M. Cannone and G. Wu in \cite{Cannone-Wu-2012} proved the global well-posedness for the Navier-Stokes system for small initial data in Fourier-Herz spaces $\dot{{\cal B}}^{-1}_{r}(\mathbb{R}^{3})$, for $1\leq r\leq 2$, corresponding to the Fourier-Besov spaces $\dot{{\cal FB}}^{-1}_{1,r}(\mathbb{R}^{3})$. In the work \cite{Iwabuchi-Takada-2014} T. Iwabuchi and R. Takada proved a global well-posedness result for the Navier-Stokes-Coriolis system for small initial data in $\dot{{\cal FB}}^{-1}_{1,2}(\mathbb{R}^{3})$, whose smallness condition does not depend on the Coriolis parameter, and also they proved the ill-posedness for initial data in $\dot{{\cal FB}}^{-1}_{1,q}(\mathbb{R}^{3})$, for any $2< q\leq \infty$. D. Fang, B. Han and M. Hieber in the work \cite{Fang-Han-Hieber-Global-Navier-Stokes-Coriolis-Fourier-Besov-2015} complement those results for small initial data in $\dot{{\cal FB}}^{2-\frac{3}{p}}_{p,r}(\mathbb{R}^{3})$, for $1<p\leq \infty$ and $1\leq r\leq\infty$.

\medskip
Related to Besov spaces, we recall the following works. T. Kato \cite{Kato-strong-1984} proved the existence of a 
unique global mild solution for small initial data in $L^{d}(\mathbb{R}^{d})$ by considering solutions spaces involving 
assymptotic properties and an auxiliary $L^{p}(\mathbb{R}^{d})$-norm, for $p\geq d$. In the work 
\cite{Planchon-NS-1996} F. Planchon proved the existence of a unique global strong solution of the Navier-Stokes equations, 
for initial data $u_0\in H^{s}(\mathbb{R}^{3})$, for $s>0$, satisfying a smallness condition in $\dot{{\cal B}}^{-\frac{1}{4}}_{4,\infty}(\mathbb{R}^{3})$. 
A similar result was obtained for initial data $u_{0}\in L^{p}(\mathbb{R}^{3})$, for $p>\frac{3}{2}$, by considering a smallness assumption 
in $\dot{{\cal B}}^{-1+\frac{3}{2p}}_{2p,\infty}(\mathbb{R}^{3})$. In the work \cite{Cannone-NS-generalization-1997} M. Cannone proved the 
existence of a unique global mild solution for initial data  $u_{0}\in L^{3}(\mathbb{R}^{3})$ satisfying a smallness assumption in the Besov space $\dot{{\cal B}}^{-1+\frac{3}{p}}_{p,\infty}(\mathbb{R}^{3})$, for $3<p\leq 6$. 
In the work \cite{Planchon-NS-1998} F. Planchon proved the existence of a unique global solution of the Navier-Stokes equations for initial 
data $u_{0}\in\dot{{\cal B}}^{-1+\frac{3}{p}}_{p,\infty}(\mathbb{R}^{3})$, with $p>3$, satisfying a smallness condition in $u_{0}\in\dot{{\cal B}}^{-1+\frac{3}{q}}_{q,\infty}(\mathbb{R}^{3})$, for some $q>p$. 
Moreover, by considering extra-conditions they show that the asymptotic behaviour of these solutions is near to self-similar solutions. 
In the work \cite{Furioli-Lemarie-Terrano-limite-2000} G. Furioli, P.G. Lemarie and E. Terrano proved the existence of a unique global mild solution for small initial 
data in the so-called limit spaces, including the class of Besov spaces $\dot{{\cal B}}^{-1+\frac{3}{p}}_{p,q}(\mathbb{R}^{3})$, for $1\leq p<3$ and 
$1\leq q<\infty$. In the work \cite{Bourgain-Pavlovic-2008} J. Bourgain and N. Pavlovi\'c proved that the solution map for the Navier-Stokes equations is 
not continuous in $\dot{{\cal B}}^{-1}_{\infty,\infty}(\mathbb{R}^{3})$ at the origin, which implies the ill-posedness of the Navier-Stokes equations for initial 
data $u_{0}\in\dot{{\cal B}}^{-1}_{\infty,\infty}(\mathbb{R}^{3})$. In the modern textbook (\cite{Bahouri-Chemin-Danchin-Fourier-PDE-2011}, 
Subsection 5.6.1) H. Bahouri, J-Y. Chemin, R. Danchin revisited these kind of results and established the global well-posedness for initial data 
$u_{0}\in\dot{{\cal B}}^{-1+\frac{3}{p}}_{p,\infty}(\mathbb{R}^{3})$, for $1\leq p<\infty$. 
Finally, we mention the work \cite{Koch-Tataru-2001} by H. Koch and D. Tataru.
They proved in that paper the global well-posedness of the Navier-Stokes system for small initial data in $BMO^{-1}(\mathbb{R}^{3})$ spaces, 
which is up to the present, one of the maximal classes for this kind of results. 

\medskip
From all of the previous cited works, one may say that, the use of Besov spaces (also Fourier-Besov spaces) has led to various advancements 
in existence, uniqueness, and regularity results for weak 
and mild solutions, particularly in spaces that allow low regularity initial data. Moreover, Besov spaces often allow us to establish 
global existence of solutions in the critical space (where the scaling invariance holds) for small initial data, which is 
a crucial aspect in fluid mechanics and turbulence theory.

\medskip
Albeit, the theory of  Besov type spaces, for instance, anisotropic Besov spaces, Besov spaces based on mixed-Lebesgue spaces, 
Besov spaces based on mixed-Lorentz spaces, etc. are under development. These spaces are extensions of classical Besov spaces,
which are designed to handle directional or spatially dependent irregularities. Indeed, in many physical and mathematical problems, 
the solution may exhibit different smoothness properties in different spatial directions. Then, mixed-norm Besov type spaces allow 
separate control over regularity along each spatial direction. This flexibility is crucial for fluid dynamics problems, 
especially when flows have directional structures, like in Navier-Stokes equations. 
In particular, we cite below two recently works in this direction.

\medskip
First, T.Phan \cite{Phan-mixed-2020} proved the global existence of a unique mild solution of the Navier-Stokes equations for small 
initial data $u_{0}\in L^{\vec{p}}(\mathbb{R}^{d})$, with $\vec{p}= (p_1, \ldots, p_d) \in [2,\infty)^{d}$, extending the previous result by Kato \cite{Kato-strong-1984} 
(for mixed-$L^{\vec{p}}(\mathbb{R}^{d})$ spaces, see \cite{Benek-Panzone-mixed-Lp-1961}). The novelty for this approach was the 
boundedness for the Riesz transform on $L^{\vec{p}}(\mathbb{R}^{d})$ and thus for the Leray operator. The solution spaces allows 
functions which decay to zero as $\mid x\mid$ tends to infinity, by considering different rates according to each spatial direction. 
Moreover, this result provides the persistence of anisotropic properties along the evolution of the system.
Second, 
in the work \cite{Khai-Tri-Sobolev-Lorentz-Navier-Stokes-2014}, D.Q. Khai and N.M. Tri proved the existence of a unique global mild solution of the Navier-Stokes system for small initial data in the Besov spaces based on mixed-Lorentz spaces $L^{\vec{q},\vec{r}}(\mathbb{R}^{d})$, or more compactly mixed-norm Besov-Lorentz spaces, that is, $u_{0}\in {\cal \dot{B}}^{-\frac{2}{p}}_{L^{\vec{q},\vec{r}},p}(\mathbb{R}^{d})$, for $\vec{r}\in [1,\infty)^{d}$, $\vec{q}\in (2,\infty)^{d}$ and $2<p<\infty$ such that $\frac{-2}{p}=-1+\displaystyle{\sum_{i=1}^{d}}\frac{1}{q_{i}}$, which implies $\displaystyle{\sum_{i=1}^{d}}\frac{1}{q_{i}}>0$. In particular, for $\vec{q}=\vec{r}$, they cover the mixed-norm Besov-Lebesgue space ${\cal \dot{B}}^{-\frac{2}{p}}_{\vec{q},p}(\mathbb{R}^{d})$, for  $\vec{q}\in (2,\infty)^{d}$ and $2<p<\infty$ such that $\frac{-2}{p}=-1+\displaystyle{\sum_{i=1}^{d}}\frac{1}{q_{i}}$.

\subsection{Main Results}
In this work we proved two results on global well-posedness for the Navier-Stokes system for small initial data in Besov type spaces based on 
mixed-Lebesgue spaces, namely, mixed-norm Besov-Lebesgue spaces (see Definition \ref{DefBesovType}) and mixed-norm Fourier-Besov-Lebesgue spaces (see Definition \ref{DefFourier-BesovType}). The former theorem proves the existence of a unique global mild solution for small initial data in the mixed-norm Besov-Lebesgue spaces 
\begin{equation*}
\dot{{\cal B}}^{\sigma}_{\vec{p},q}(\mathbb{R}^{d}),\,\,\mbox{where}\,\,\sigma=-1+\displaystyle{\sum_{i=1}^{d}}\frac{1}{p_{i}},\,\, \mbox{for}\,\, 1\leq q\leq \infty\,\, \mbox{and}\,\,\vec{p}\in [1,\infty]^{d}\,\,\mbox{such that}\,\,\displaystyle{\sum_{i=1}^{d}}\frac{1}{p_{i}}>0.
\end{equation*}
In particular, we need at least some $p_{i}\in [1,\infty)$. In fact, for $p_{1}=\dots=p_{d}=\infty$, 
the system \eqref{BQ} is ill-posed, see \cite{Bourgain-Pavlovic-2008}. Therefore, we have extended the results by 
G. Furioli, P.G. Lemarie and E. Terrano \cite{Furioli-Lemarie-Terrano-limite-2000}, D.G. Khai and N.M. Tri \cite{Khai-Tri-Sobolev-Lorentz-Navier-Stokes-2014} 
and by H. Bahouri, J-Y. Chemin and R. Danchin (\cite{Bahouri-Chemin-Danchin-Fourier-PDE-2011}, Subsection 5.6.1). 
It is also important to observe that 
 \begin{equation*}
 \begin{split}
&\dot{{\cal B}}^{-1+\frac{3}{q}}_{q,\infty}(\mathbb{R}^{3})\hookrightarrow BMO^{-1}(\mathbb{R}^{3})\,\,\mbox{for}\, 3<q<\infty\,\,\mbox{and}\\
&\dot{{\cal B}}^{-1+\frac{3}{q}}_{q,\infty}(\mathbb{R}^{3})\hookrightarrow \dot{{\cal B}}^{\sigma}_{\vec{p},\infty}(\mathbb{R}^{3}),\,\,\mbox{where}\,\,\sigma=-1+\displaystyle{\sum_{i=1}^{3}}\frac{1}{p_{i}},\,\,\mbox{for}\,\,1\leq q<\infty,\,\, (q,q, q)\leq \vec{p}\leq (\infty,\infty,\infty).    
\end{split}
 \end{equation*}
Then, for the considered parameters 
$$
   \text{$3<q<\infty$ and $(q,q, q)\leq \vec{p}\leq (\infty,\infty,\infty)$, 
such that $\vec{p}\neq (\infty,\infty,\infty)$},
$$
there is no clear inclusion between $BMO^{-1}(\mathbb{R}^{3})$ 
and $\dot{{\cal B}}^{\sigma}_{\vec{p},\infty}(\mathbb{R}^{3})$, for $\sigma=-1+\displaystyle{\sum_{i=1}^{3}}\frac{1}{p_{i}}$, 
and thus is expected that our result provides a new class of initial data for the global well-posedness of the Navier-Stokes equations, 
which are not included in $BMO^{-1}(\mathbb{R}^{3})$, however included in $\dot{{\cal B}}^{-1}_{\infty,\infty}(\mathbb{R}^{3})$. 
Moreover, in the recent book (\cite{Lemarie-NS-century-2018}, Lemma 7.10) Lemari\'e-Rieusset established the equivalence
\begin{equation*}
f\in \dot{{\cal B}}^{s}_{E,\infty}\iff 
\displaystyle{\sup_{t>0}}\,t^{\frac{\mid s\mid}{2}}\parallel S(t)f\parallel_{E}<\infty,
\end{equation*}
where $s<0$ and $E$ is a Banach space which is a subset of the set of tempered distributions $S'(\mathbb{R}^{d})$ and $E$ is stable under convolution with $L^{1}(\mathbb{R}^{d})$, that is,
\begin{equation*}
\parallel f\ast g\parallel_{E}\leq \parallel f\parallel_{1}\parallel g\parallel_{E},\,\,\mbox{for all}\,\, f\in L^{1}(\mathbb{R}^{d})\,\,\mbox{and}\,\,g\in E.
\end{equation*}
In particular the space $E=L^{\vec{p}}(\mathbb{R}^{d})$, with $\vec{p}\in [1,\infty]^{d}$ satisfies these conditions. 
This lemma allows to see that, we also complement the previous and recent global well-posedness result in \cite{Phan-mixed-2020}.

\medskip
Moreover, as a direct implication we can say that, for the next family of Besov type spaces $\dot{{\cal B}}^{-1+\frac{1}{p}}_{(\infty,\infty,p),\infty}$, for $1\leq p\leq \infty$, 
our result provides a new optimal class for the well-posedness for $1\leq p< \infty$, since for $p=\infty$ the system is ill-posedness as proved in 
the well-known result \cite{Bourgain-Pavlovic-2008} by J. Bourgain and N. Pavlovi\'c. 

\medskip
In the second main theorem of this paper, we proved the existence of a unique global mild solution for the Navier-Stokes equations, for initial data in the mixed-norm Fourier-Besov-Lebesgue spaces $\dot{{\cal FB}}^{s}_{\vec{p},q}(\mathbb{R}^{d})$,
\begin{equation*}
\mbox{where}\,\,s=-1+\displaystyle{\sum_{i=1}^{d}}\left(1-\frac{1}{p_{i}}\right),\,\, \mbox{for}\,\, 1\leq q\leq \infty\,\, \mbox{and}\,\,\vec{p}\in [1,\infty]^{d}\,\,\mbox{such that}\,\,\displaystyle{\sum_{i=1}^{d}}\left(1-\frac{1}{p_{i}}\right)>0,    
\end{equation*}
 which allow to partly extend previous results by P. Konieczny and T. Yoneda \cite{Konieczny-Yoneda-2011}, T. Iwabuchi and R. Takada \cite{Iwabuchi-Takada-2014}, Z. Lei and F.H. Lin \cite{Lein-Lin-2011} and D. Fang, B. Han and M. Hieber \cite{Fang-Han-Hieber-Global-Navier-Stokes-Coriolis-Fourier-Besov-2015}.  
 
 Finally, as a particular case we have that, for the family of Fourier-Besov type spaces $\dot{{\cal FB}}^{-\frac{1}{p}}_{(1,1,p),\infty}$, where $1\leq p\leq \infty$, 
 our result provide a new optimal class for the well-posedness for $1< p\leq \infty$, 
 since for $p=1$ the system \eqref{BQ}
 is ill-posed due to the well-known result by T. Iwabuchi and R. Takada \cite{Iwabuchi-Takada-2014}.


\bigskip
The structure of this paper is the following:

\medskip
In Section \ref{preliminar}, we first define mixed-norm Besov-Lebesgue spaces, mixed-norm Fourier-Besov-Lebesgue spaces and we 
established some primordial properties of these spaces. We also recall the definitions of Stokes semigroup 
and mild solutions for the Navier-Stokes system. Then, in Section \ref{main-theorems}, we just set our main 
results on global well-posedness in the context of mixed-norm Besov-Lebesgue type spaces.
In Section \ref{main-tools}, we recall two of the main technical tools for this approach, 
as being the Bony's paraproduct and an abstract fixed point scheme. In this section we also established the linear estimates.
In Section \ref{bilinear-estimates-s}, we make use of linear estimates and Bony's paraproduct in order to get the bilinear estimates.
In the last section, we indicated the procedure in order to apply the fixed point scheme.

\section{Preliminaries}\label{preliminar}

In this section, we fix the notations and collect some preliminary results. 
First, let $\Omega \subset \mathbb{R}^d$ be an open set. 
We denote by $dx, d\xi$, etc. 
the Lebesgue measure on $\Omega$ and by $L^p(\Omega)$, $p \in [1,+\infty)$, the set of $p$-summable functions 
with respect to the Lebesgue measure, also $L^\infty(\Omega)$ is 
the set of measurable functions in which its absolute value has the essential supremum finite. 
As is customary, the symbols $\mathcal{S}(\mathbb{R}^d)$ denotes the Schwartz class, and
$\mathcal{S}^{\prime}(\mathbb{R}^d)$ denotes the set of tempered distributions. We denote by $\mathcal{F} \varphi(\xi) \equiv \widehat{\varphi}(\xi)\equiv\varphi^{\wedge}(\xi)$
the Fourier Transform of $\varphi$, which is an
isometry in $L^2(\mathbb{R}^d)$.

\medskip
Let ${\cal A}=\{\eta\in \mathbb{R}^{d};\frac{3}{4}\leq \mid \eta\mid\leq \frac{8}{3}\}$ 
be a given annulus and $\phi\in {\cal S}(\mathbb{R}^{d})$ be a radially symmetric function for which $\supp \hat{\phi}\subset {\cal A}$, satisfying the identity
\begin{equation*}
\displaystyle{\sum_{l\in\mathbb{Z}}}\phi(2^{-l}\eta)=1,\,\,\mbox{for each}\,\,\eta\neq 0.    
\end{equation*}

Let us denote $\phi_l(\eta)=\phi(2^{-l}\eta)$ and $\zeta_l(\eta)=\displaystyle{\sum_{l'\leq l-1}\phi_{l'}(\eta)}$. Now, we recall the standard localization operators, defined by 
\begin{equation*}
\dot{\Delta}_{l}g=\phi_l(D)g,\,\,\,\dot{S}_lg=\displaystyle{\sum_{l'\leq l-1}}\dot{\Delta}_{l'}g=\zeta_{l}(D)g,\,\,\mbox{where}\,\,l\in\mathbb{Z},    
\end{equation*}
where $h(D)$ denotes the pseudo-differential operator defined by $h(D)u={\cal F}^{-1}[h(\eta)\hat{u}]$. Here ${\cal F}^{-1}$ denotes the inverse Fourier transform.

The main properties of these operators are the next almost orthogonal ones, that is to say
\begin{equation*}
\begin{split}
\dot{\Delta}_l \dot{\Delta}_{l'}g&=0,\,\,\mbox{if}\,\,\mid l-l'\mid\geq 2,\,\,\mbox{and}\\
\dot{\Delta}_l\left(\dot{S}_{l'-1}f\dot{\Delta}_{l'}g \right)&=0, \,\,\mbox{if}\,\,\mid l-l'\mid\geq 5.
\end{split}
\end{equation*}

\subsection{Mixed-norm Besov-Lebesgue spaces}
In this section we recall the definition of mixed-norm Besov-Lebesgue spaces and its main properties as Bernstein-type lemmas, embedding inclusions and completeness. 

We first recall the mixed-Lebesgue spaces and its main properties.
 \begin{defi}(Mixed-Lebesgue spaces)
 Let us consider $\vec{p}=(p_{1},\dots, p_{d})\in [1,\infty]^{d}$. The Mixed-Lebesgue spaces $L^{\vec{p}}(\mathbb{R}^{d})$ is the set of measurable functions $f$ defined on $\mathbb{R}^{d}$ such that 
 \begin{equation*}
   \parallel f\parallel_{L^{\vec{p}}(\mathbb{R}^{d})}=\left(
   \displaystyle{\int_{\mathbb{R}}}\left( \displaystyle{\int_{\mathbb{R}}}\dots \left(\displaystyle{\int_{\mathbb{R}}}\mid f(x_{1},\dots, x_{d})\mid^{p_{1}}dx_{1}\right)^{\frac{1}{p_{1}}}\dots dx_{d-1}\right)^{\frac{p_{d}}{p_{d-1}}} dx_{d}\right)^{\frac{1}{p_{d}}}\,\,\mbox{is finite}.
 \end{equation*}
 Here we consider integrals if $\vec{p}\in [1,\infty)^{d}$ and the $i-$th essential sup norm if some $p_{i}=\infty$. 
 \end{defi}
Next lemma summarize some main properties on mixed-Lebesgue spaces.
\begin{lema} 
 \item [(i.)] For $\vec{p}_{1},\vec{p}_{2},\vec{p}_{3}\in [1,\infty]^{d}$, such that $\frac{1}{\vec{p}_{3}}=\frac{1}{\vec{p}_{1}}+\frac{1}{\vec{p}_{2}}$ we have
\begin{equation*}
 \parallel g_1\cdot g_2\parallel_{L^{\vec{p}_{3}}(\mathbb{R}^{d})}\leq    \parallel g_1\parallel_{L^{\vec{p}_{1}}(\mathbb{R}^{d})}\cdot
 \parallel g_2\parallel_{L^{\vec{p}_{2}}(\mathbb{R}^{d})},
\end{equation*}
for all $g_{1}\in L^{\vec{p_{1}}}(\mathbb{R}^{d})$ and $g_{2}\in L^{\vec{p_{2}}}(\mathbb{R}^{d})$.
\item  [(ii.)] For $\vec{p}\in[1,\infty]^{d}$ we have 
\begin{equation*}
 \parallel \phi\ast g\parallel_{L^{\vec{p}}(\mathbb{R}^{d})}\leq \parallel \phi\parallel_{L^{1}(\mathbb{R}^{d})}\parallel g\parallel_{L^{\vec{p}}(\mathbb{R}^{d})},
 \end{equation*}
 for all $\phi \in L^{1}(\mathbb{R}^{d})$ and $g\in L^{\vec{p}}(\mathbb{R}^{d})$.
 \item [(iii.)] The normed space $\left(L^{\vec{p}}(\mathbb{R}^{d}),\parallel \cdot\parallel_{L^{\vec{p}}(\mathbb{R}^{d})}\right)$ is complete.
\end{lema}
Here, we denote $\frac{1}{\vec{p}}=(\frac{1}{p_{1}},\dots,\frac{1}{p_{d}})$, where $\vec{p}=(p_{1},\dots,p_{d})\in[1,\infty]^{d}$.
For further properties on these spaces see the work \cite{Benek-Panzone-mixed-Lp-1961}, by A. Benedek and R. Panzone. 
 
\bigskip
In the modern definition, mixed-norm Besov spaces were introduced recently in the work \cite{Cleanthous-Georgiadis-Nielsen-mixed-norm-2016}, 
where they established some on their main properties as well as the discrete decomposition of homogeneous mixed-norm Besov spaces. 
To see another related studies we address the reader to \cite{Cleanthous-Georgiadis-Nielsen-discrete-2017} and \cite{Georgiadis-Nielsen-Pseudo-2016}, also 
for usual Besov spaces see \cite{Sawano-Besov-2018} and \cite{Bahouri-Chemin-Danchin-Fourier-PDE-2011}.
Before we introduce the mixed-norm Besov-Lebesgue spaces, similar in the spirit, however slightly different to the previous mixed-norm Besov spaces, 
we set 
$$
   S_{h}'= \Big \{g\in S'(\mathbb{R}^{d});\displaystyle{\lim_{\lambda\rightarrow \infty}}\parallel \theta(\lambda D)g\parallel_{\infty}=0,\,\,\mbox{for any}\,\,\theta\in C^{\infty}_{0}(\mathbb{R}^{d}) \Big\}.
 $$  
%
\begin{defi}(Mixed-norm Besov-Lebesgue spaces)
\label{DefBesovType} 
Let us consider three parameters 
$$
    \text{$\vec{p}=(p_1,\dots,p_d)\in [1,\infty]^{d}$, $q\in [1,\infty]$ and $\sigma\in\mathbb{R}$.}
$$    
The $d$-dimensional homogeneous mixed-norm Besov-Lebesgue spaces are defined as
\begin{itemize}
\item [(i.)] For $\vec{p}\in [1,\infty]^{d}$, $q\in[1,\infty)$ and $\sigma\in\mathbb{R}$,
\begin{equation}
\dot{{\cal B}}^{\sigma}_{\vec{p},q}(\mathbb{R}^{d})=\lbrace
g\in {\cal S}_{h}';\parallel g\parallel_{\dot{{\cal B}}^{\sigma}_{\vec{p},q}(\mathbb{R}^{d})}=
\left(\displaystyle{\sum_{l\in\mathbb{Z}}2^{l\sigma q}\parallel \dot{\Delta}_{l}g\parallel^{q}_{L^{\vec{p}}(\mathbb{R}^{d})}}
\right)^{\frac{1}{q}}\,\,\mbox{is finite}
\rbrace.
\end{equation}
\item [(ii.)]  For $\vec{p}\in [1,\infty]^{d}$, $q=\infty$ and $\sigma\in\mathbb{R}$,
\begin{equation}
\dot{{\cal B}}^{\sigma}_{\vec{p},\infty}(\mathbb{R}^{d})=\lbrace
g\in {\cal S}_{h}';\parallel g\parallel_{\dot{{\cal B}}^{\sigma}_{\vec{p},\infty}(\mathbb{R}^{d})}=\displaystyle{\sup_{l\in\mathbb{Z}} 2^{l\sigma}\parallel \dot{\Delta}_{l}g\parallel_{L^{\vec{p}}(\mathbb{R}^{d})}}\,\,\mbox{is finite}
\rbrace.
\end{equation}
\end{itemize}
We also called these spaces as the $d$-dimensional homogeneous Besov mixed-Lebesgue spaces or the $d$-dimensional homogeneous Besov spaces based on mixed-Lebesgue spaces.  
\end{defi}
\begin{lema}{(Main properties related to mixed-norm Besov-Lebesgue spaces)}\label{properties-Besov-spaces-mixed}
\begin{itemize}
\item [(i.)]  Let ${\cal A}=\{\eta\in \mathbb{R}^{d}; r_{1}\leq \mid\eta\mid\leq r_{2}\}$, for some $0<r_{1}<r_{2}$, be a given annulus and 
$$
\text{$B=\{\eta\in \mathbb{R}^{d}; \mid \eta\mid\leq r\}$, for some $r>0$, be a closed ball.} 
$$
Then, there exist a positive constant $C$, such that, for any nonnegative integer $k$ and a couple 
$\vec{p},\vec{q}\in\left([1,\infty]^{d}\right)^{2}$, with $\vec{1}\leq\vec{p}\leq \vec{q}$, it follows for any function $u\in L^{\vec{p}}(\mathbb{R}^{d})$ that:
\begin{equation*}
 \begin{split}
   &(i.1.)\,\,\supp \hat{u}\subset \lambda B\implies\parallel D^{k}u\parallel_{L^{\vec{q}}(\mathbb{R}^{d})}=\displaystyle{\sup_{\mid \alpha\mid=k}}\parallel \partial^{\alpha}u\parallel_{L^{\vec{q}}(\mathbb{R}^{d})}\leq C^{k}\lambda^{k_{0}}\parallel u\parallel_{L^{\vec{p}}(\mathbb{R}^{d})}.\\
   &(i.2.)\,\, \supp\hat{u}\subset \lambda {\cal A}\implies C^{-k-1}\lambda^{k}\parallel u\parallel_{L^{\vec{p}}(\mathbb{R}^{d})}\leq \parallel D^{k}u\parallel_{L^{\vec
   p}(\mathbb{R}^{d})}\leq C^{k+1}\lambda^{k}\parallel u\parallel_{L^{\vec{p}}(\mathbb{R}^{d})},  
 \end{split}   
\end{equation*}
where $k_{0}=k+\displaystyle{\sum_{i=1}^{d}}\left(
   \frac{1}{p_i}-\frac{1}{q_i}\right)$.
\item[(ii.)] Let us consider $\vec{p},\vec{q}\in[1,\infty]^{d}$ such that $\vec{1}\leq\vec{p}\leq\vec{q}$ and $1\leq a_{1}\leq a_{2}\leq \infty$. For $\sigma_{2}<\sigma_{1}$ 
with $\sigma_{1}+\displaystyle{\sum_{i=1}^{d}}\frac{1}{q_{i}}=\sigma_{2}+\displaystyle{\sum_{i=1}^{d}}\frac{1}{p_{i}}$, we have
\begin{equation*}
 \dot{{\cal B}}^{\sigma_{1}}_{\vec{p},a_{1}}\hookrightarrow
 \dot{{\cal B}}^{\sigma_{2}}_{\vec{q},a_{2}},
\end{equation*}
that is, the above inclusion is continuous. Here, $\vec{1}=(1,\dots,1)$ and $\vec{p}\leq\vec{q}$ means $p_{1}\leq q_{1},\dots, p_{d}\leq q_{d}$ if $\vec{p}=(p_{1},\dots,p_{d})$ and $\vec{q}=(q_{1},\dots,q_{d})$.
\item [(iii.)] For $\vec{p}\in[1,\infty]^{d}$, if $\sigma<\displaystyle{\sum_{i=1}^{d}}\frac{1}{p_{i}}$ and $1\leq q\leq \infty$ or $\sigma=\displaystyle{\sum_{i=1}^{d}}\frac{1}{p_{i}}$ and $q=1$, then $\dot{{\cal B}}^{\sigma}_{\vec{p},q}(\mathbb{R}^{d})$ is a Banach space endowed with the norm $\parallel \cdot\parallel_{\dot{{\cal B}}^{\sigma}_{\vec{p},q}(\mathbb{R}^{d})}$.
\end{itemize}   
\end{lema}
\subsection{Mixed-norm Fourier-Besov-Lebesgue spaces}
Fourier-Besov spaces were introduced in the work of P. Konieczny and T. Yoneda \cite{Konieczny-Yoneda-2011} to study the dispersive effect of the Coriolis force for both the stationary and the non-stationary Navier-Stokes system. To the best of our knowledge the next definition of Fourier-Besov type spaces is new in the context of mixed-Lebesgue spaces.

\begin{defi}(Mixed-norm Fourier-Besov-Lebesgue spaces)\label{DefFourier-BesovType}
Let us consider three parameters 
$$
    \text{$\vec{p}=(p_1,\dots,p_d)\in [1,\infty]^{d}$, $q\in [1,\infty]$ and $s\in\mathbb{R}$.}
$$ 
The $d$-dimensional homogeneous mixed-norm Fourier-Besov-Lebesgue spaces are defined as
\begin{itemize}
\item [(i.)] For $\vec{p}\in [1,\infty]^{d}$, $q\in[1,\infty)$ and $s\in\mathbb{R}$,
\begin{equation}
\dot{{\cal FB}}^{s}_{\vec{p},q}(\mathbb{R}^{d})=\lbrace
g\in {\cal S}_{h}';\parallel g\parallel_{\dot{{\cal FB}}^{s}_{\vec{p},q}(\mathbb{R}^{d})}=
\left(\displaystyle{\sum_{l\in\mathbb{Z}}2^{lsq}\parallel \phi_{l}\hat{g}\parallel^{q}_{L^{\vec{p}}(\mathbb{R}^{d})}}
\right)^{\frac{1}{q}}\,\,\mbox{is finite}
\rbrace.
\end{equation}
\item [(ii.)]  For $\vec{p}\in [1,\infty]^{d}$, $q=\infty$ and $s\in\mathbb{R}$,
\begin{equation}
\dot{{\cal FB}}^{s}_{\vec{p},\infty}(\mathbb{R}^{d})=\lbrace
g\in {\cal S}_{h}';\parallel g\parallel_{\dot{{\cal FB}}^{s}_{\vec{p},\infty}(\mathbb{R}^{d})}=\displaystyle{\sup_{l\in\mathbb{Z}} 2^{ls}\parallel \phi_{l}g\parallel_{L^{\vec{p}}(\mathbb{R}^{d})}}\,\,\mbox{is finite}
\rbrace. 
\end{equation}
\end{itemize}
We also called these spaces as the $d$-dimensional homogeneous Fourier-Besov mixed-Lebesgue spaces or the $d$-dimensional homogeneous Fourier-Besov spaces based on mixed-Lebesgue spaces.
\end{defi}
\begin{lema}{(Main properties related to mixed-norm Fourier-Besov-Lebesgue spaces)}\label{properties-Fourier-Besov-mixed}
\begin{itemize}
\item [(i.)] Let $\vec{p},\vec{q}\in[1,\infty]^{d}$ such that $\vec{p}\leq \vec{q}$ and let $g$ be any map with $\hat{g}\in L^{\vec{q}}(\mathbb{R}^{d})$ such that 

\begin{equation*}
 \supp \hat{g}\subset\{
 \xi_{1}\in\mathbb{R};\mid \xi_{1}\mid \leq A_{1}2^{i_1}\}\times\cdots\times\{
 \xi_{d}\in\mathbb{R};\mid \xi_{d}\mid \leq A_{d}2^{i_d}\},   
\end{equation*}
for some real constants $A_1,\cdots, A_d$ and integers $i_{1},\cdots, i_{d}\in \mathbb{Z}$. Then,
\begin{equation*}
 \parallel\left(\xi_1\right)^{\beta_1}\cdots \left(\xi_d\right)^{\beta_d}\hat{g}\parallel_{L^{\vec{p}}(\mathbb{R}^{d})}\leq C2^{k_{0}}  \parallel\hat{g}\parallel_{L^{\vec{q}}(\mathbb{R}^{d})},
\end{equation*}
where $k_{0}=\displaystyle{\sum_{k=1}^{d}}i_{k}\left(\mid \beta_{k}\mid+\frac{1}{p_{k}}-\frac{1}{q_{k}}\right)$ and $C$ is a constant depending on $A_{1},\cdots,A_{d}$, for any indexes $\beta_{1},\cdots,\beta_{d}$.
\item[(ii.)] Let us consider $\vec{p},\vec{q}\in[1,\infty]^{d}$, $\vec{q}\leq\vec{p}$ and $1\leq a_{1}\leq a_{2}\leq \infty$. If we consider $s_{1},s_{2}\in\mathbb{R}$ such that $s_{1}\geq s_{2}$ and $s_{1}+\displaystyle{\sum_{i=1}^{d}}\frac{1}{p_{i}}=s_{2}+\displaystyle{\sum_{i=1}^{d}}\frac{1}{q_{i}}$, then we have
\begin{equation*}
 \dot{{\cal FB}}^{s_{1}}_{\vec{p},a_{1}}\hookrightarrow
 \dot{{\cal FB}}^{s_{2}}_{\vec{q},a_{2}},
\end{equation*}
that is, the above inclusion is continuous.
\item [(iii.)] For $\vec{p}\in[1,\infty]^{d}$, if $s<\displaystyle{\sum_{i=1}^{d}}\left( 1-\frac{1}{p_{i}}\right)$ and $1\leq q\leq \infty$ or $s=\displaystyle{\sum_{i=1}^{d}}\left( 1-\frac{1}{p_{i}}\right)$ and $q=1$ then $\dot{{\cal FB}}^{s}_{\vec{p},q}(\mathbb{R}^{d})$ is a Banach space endowed with the norm $\parallel \cdot\parallel_{\dot{{\cal FB}}^{s}_{\vec{p},q}(\mathbb{R}^{d})}$.
\end{itemize}   
\end{lema}
\subsection{The Stokes semigroup and Mild solutions}
We shall denote by $S_{\nu}(t)=e^{\nu t\Delta}$, for $t>0$, the Stokes semigroup whose Fourier representation is given by 
\begin{equation*}
    [S_{\nu}(t)g]^{\wedge}(\eta)=e^{-\nu t\mid \eta\mid^{2}}\hat{g}(\eta).  
\end{equation*}
Let us denote by $\mathbb{P}$ the Leray projector, which in Fourier variables is given by its symbols $\hat{\mathbb{P}}$, that is, 
\begin{equation*}
 \left(\hat{\mathbb{P}}(\eta)\right)_{i,j}=\delta_{i,j}-\frac{\eta_{i}\eta_{j}}{\mid \eta\mid^{2}}, \,\,\mbox{for each}\,\, i,j=1,\dots, d,   
\end{equation*}
where $\delta_{i,j}$ denotes the Kronecker delta.

The mild formulation for the Navier-Stokes equations is given by the next abstract equation
\begin{equation}\label{mild-formulation}
u(t)=S_{\nu}(t)u_{0}+B(u,u)(t),    
\end{equation}
for all $t\in[0,\infty)$, where the bilinear operator $B(\cdot,\cdot)$ is given by 
\begin{equation*}
 B(v,w)(t)=\displaystyle{\int_{0}^{t}}S_{\nu}(t-t')\mathbb{P}\mbox{div}\left(v\otimes w\right)(t')\,dt',   
\end{equation*}
which is given in Fourier variables by 
\begin{equation*}
[B(v,w)]^{\wedge}(\eta,t)=\displaystyle{\int_{0}^{t}}e^{-\nu(t-t')\mid\eta\mid^{2}}\hat{\mathbb{P}}(\eta)[i\eta\cdot \left(v\otimes w\right)^{\wedge}(\eta,t')]\,dt', 
\end{equation*}
where the components of the matrix $(u\otimes v)^{\wedge}$ and the vector $i\eta\cdot\left(v\otimes w\right)^{\wedge}$ are given by 
\begin{equation*}
  \left(v\otimes w\right)^{\wedge}_{i,j}=\hat{v}_{i}\ast\hat{w}_{j}, \,\,\mbox{for}\,\,i,j=1,\dots,d,\,\,\mbox{and}\,\,\left(i\eta\cdot\left(v\otimes w\right)^{\wedge}\right)_{l}=i\eta\cdot\mbox{row}_{l}\,\left(v\otimes w\right)^{\wedge},\,\,\mbox{for}\,\,l=1,\dots,d.  
\end{equation*}
\begin{defi}
 Let $Z$ be a suitable Banach space. We say that $u\in Z$ is a mild solution for the Navier-Stokes equations if $u$ satisfies the mild formulation (\ref{mild-formulation}) above.
 \begin{obs}
In the context of the main theorems below $Z={\cal L}^{\infty}\left(I;{\dot{{\cal B}}^{\sigma}_{\vec{p},q}(\mathbb{R}^{d})}\right)\bigcap {\cal L}^{1}\left(I;{\dot{{\cal B}}^{\sigma+2}_{\vec{p},q}(\mathbb{R}^{d})}\right)$, so that in order to get $S_{\nu}(\cdot)u_{0}\in Z$, we required the initial data $u_{0}\in \dot{{\cal B}}^{\sigma}_{\vec{p},q}(\mathbb{R}^{d})$. Similar remark for Fourier-Besov mixed-Lebesgue spaces.  
 \end{obs}
\end{defi}
Let us recall the space ${\cal L}^{a}(I;\dot{{\cal FB}}^{s}_{\vec{p},q}(\mathbb{R}^{d}))$ of Bochner measurable functions defined on interval $I$ toward $\dot{{\cal FB}}^{s}_{\vec{p},q}(\mathbb{R}^{d})$ provided by the norm
\begin{equation*}
\parallel g\parallel_{{\cal L}^{a}(I;\dot{{\cal FB}}^{s}_{\vec{p},q}(\mathbb{R}^{d}))}=\displaystyle{\parallel
2^{js} \parallel\phi_{l}\hat{g}\parallel_{L^{a}(I;L^{\vec{p}}(\mathbb{R}^{d}))}\parallel_{l^{q}(\mathbb{Z})}}.    
\end{equation*}
A similar definition should be considered if we consider $\dot{{\cal B}}^{\sigma}_{\vec{p},q}(\mathbb{R}^{d})$ instead of $\dot{{\cal FB}}^{s}_{\vec{p},q}(\mathbb{R}^{d})$.

\section{The main theorem}
\label{main-theorems}

In this section, we state our global well-posedness results for the Navier-Stokes system. 
\begin{teo}
Let $\vec{p}\in[1,\infty]^{d}$ be a parameter such that $\displaystyle{\sum_{i=1}^{d}}\frac{1}{p_i}>0$, $\sigma=-1+\displaystyle{\sum_{i=1}^{d}}\frac{1}{p_i}$ and \newline 
$Z={\cal L}^{\infty}\left(I;{\dot{{\cal B}}^{\sigma}_{\vec{p},q}(\mathbb{R}^{d})}\right)\bigcap {\cal L}^{1}\left(I;{\dot{{\cal B}}^{\sigma+2}_{\vec{p},q}(\mathbb{R}^{d})}\right)$, where $I=(0,\infty)$. Then,
\begin{itemize}
\item [(i.)] There exist $\delta>0$ and $D_0$ such that, for small enough initial data $\parallel u_{0}\parallel_{\dot{{\cal B}}^{\sigma}_{\vec{p},q}(\mathbb{R}^{d})}\leq \delta$ there is a unique global mild solution $u$ in the ball $\{u\in Z;\parallel u\parallel_{Z}\leq \delta D_0\}$. Moreover, the map solution depends continuously on the initial data. 
\item [(ii.)] The obtained solution is weakly continuous from $[0,\infty)$ to $S'(\mathbb{R}^{d})$. 
\end{itemize}
\end{teo}
As mentioned in the introduction this theorem provides a complement and partial extension to the results by G. Furioli, P.G. Lemarie and E. Terrano \cite{Furioli-Lemarie-Terrano-limite-2000}, D.G. Khai and N.M. Tri \cite{Khai-Tri-Sobolev-Lorentz-Navier-Stokes-2014} and by H. Bahouri, J-Y. Chemin and R. Danchin (\cite{Bahouri-Chemin-Danchin-Fourier-PDE-2011}, Subsection 5.6.1).
\begin{teo}
Let $\vec{p}\in[1,\infty]^{d}$ be a parameter such that $\displaystyle{\sum_{i=1}^{d}}\left(1-\frac{1}{p_i}\right)>0$, $s=-1+\displaystyle{\sum_{i=1}^{d}}\left(1-\frac{1}{p_i}\right)$ and $X={\cal L}^{\infty}\left(I;{\dot{{\cal FB}}^{s}_{\vec{p},q}(\mathbb{R}^{d})}\right)\bigcap {\cal L}^{1}\left(I;{\dot{{\cal FB}}^{s+2}_{\vec{p},q}(\mathbb{R}^{d})}\right)$, where $I=(0,\infty)$. Then,
\begin{itemize}
\item [(i.)] There exist $\delta>0$ and $D_0$ such that, for small enough initial data $\parallel u_{0}\parallel_{\dot{{\cal FB}}^{s}_{\vec{p},q}(\mathbb{R}^{d})}\leq \delta$ there is a unique global mild solution $u$ in the ball $\{u\in X;\parallel u\parallel_{X}\leq \delta D_0\}$. Moreover, the map solution depends continuously on the initial data. 
\item [(ii.)] The obtained solution is weakly continuous from $[0,\infty)$ to $S'(\mathbb{R}^{d})$. 
\end{itemize}
\end{teo}
As observed in the introduction this theorem also provides a complement and partial extension to the results by P. Konieczny and T. Yoneda \cite{Konieczny-Yoneda-2011}, T. Iwabuchi and R. Takada \cite{Iwabuchi-Takada-2014}, Z. Lei and F.H. Lin \cite{Lein-Lin-2011} and D. Fang, B. Han and M. Hieber \cite{Fang-Han-Hieber-Global-Navier-Stokes-Coriolis-Fourier-Besov-2015}.  
\begin{obs}
In order to complement the context of these main theorems we also mention the works on anisotropic Besov type spaces, by considering initial data in ${\cal B}^{0,\frac{1}{2}}(\mathbb{R}^{3})$ \cite{Paicu-anisotrope-2005} by Paicu M., or in ${\cal B}^{-\frac{1}{2},\frac{1}{2}}_{4}(\mathbb{R}^{3})$ \cite{Paicu-Zhang-global-2011} by Paicu, M. and Zhang, P., and in Besov type spaces based on Morrey spaces, namely, in Besov-Morrey spaces ${\cal N}^{-1+\frac{d}{p}}_{p,q,\infty}(\mathbb{R}^{d})$, for $1\leq q\leq p<\infty$ and $d< p$, \cite{Kozono-Yamazaki-semilinear-1994} by Kozono, H. and Yamazaki, M. and in Fourier-Besov-Morrey spaces $\dot{{\cal FN}}^{2-\frac{3-\mu}{q}}_{q,\mu,\infty}(\mathbb{R}^{3})$, for $0\leq \mu<3$ and $1\leq q<\infty$ ($\mu\neq 0$ if $q=1$), \cite{Almeida-Ferreira-Lima-uniform-2017} by Almeida, M.F.d., Ferreira, L.C.F., and Lima, L.S.M.. In fact, it is also important to consider analogous Besov mixed-Morrey spaces and Fourier-Besov mixed-Morrey spaces.
\end{obs}
\begin{obs}
In this remark we establish some relation between both the above given theorems. We first recall the Hausdorff-Young inequality (\cite{Folland-real-1999}, p. 253),
\begin{equation*}
 \parallel f\parallel_{L^{p'}(\mathbb{R}^{d})}\leq \parallel f\parallel_{L^{p}(\mathbb{R}^{d})},\,\,\mbox{for}\,\,1\leq p\leq 2,   
\end{equation*}
where $p'$ denote the conjugate exponent to $p$. We also have, by applying repeatdly this last inequality and the Minkowski's inequality, the estimate
\begin{equation*}
 \parallel f\parallel_{L^{\vec{p}'}(\mathbb{R}^{d})}\leq \parallel f\parallel_{L^{\vec{p}}(\mathbb{R}^{d})},\,\,\mbox{for}\,\,\vec{p}=(p_{1},\dots,p_{d})\,\,\mbox{such that}\,\,1\leq p_{d}\leq \dots\leq p_{1}\leq 2,   
\end{equation*}
which implies 
\begin{equation*}
 \dot{{\cal B}}^{\sigma}_{\vec{p},q}(\mathbb{R}^{d})\hookrightarrow \dot{{\cal FB}}^{\sigma}_{\vec{p}',q}(\mathbb{R}^{d})\,\,\mbox{and}\,\,\dot{{\cal FB}}^{s}_{\vec{p},q}(\mathbb{R}^{d})\hookrightarrow \dot{{\cal B}}^{s}_{\vec{p}',q}(\mathbb{R}^{d}).    
\end{equation*}
This inclusion relationships suggest that the obtained global existence results could be the maximal in comparison between both the classes of mixed-norm Besov-Lebesgue types spaces and thus we note the same relevance for each space. For instance, we can take small initial data $u_{0}\in \dot{{\cal FB}}^{0}_{(2,1,2),q}(\mathbb{R}^{3})$, and small initial data $u_{0}\in \dot{{\cal B}}^{0}_{(2,\infty,2),q}(\mathbb{R}^{3})$, for $1\leq q\leq \infty$, in order to obtain the global well-posedness. However we can not assure
\begin{equation}
\dot{{\cal FB}}^{0}_{(2,1,2),q}(\mathbb{R}^{3})\hookrightarrow \dot{{\cal B}}^{0}_{(2,\infty,2),q}(\mathbb{R}^{3})
\end{equation}
and therefore the class $\dot{{\cal FB}}^{0}_{(2,1,2),q}(\mathbb{R}^{3})$ could be maximal between the both last classes.
\end{obs}

\section{Some auxiliary tools and linear estimates}
\label{main-tools}

In this section we recall two main analytical tools which we shall use along this work, as being the Bony's 
paraproduct and a fixed point scheme. Moreover we establish linear estimates involving the Stokes semigroup, 
which we shall use for computations in relation to the initial data and also for the auxiliary linear operator related to the bilinear term.

\subsection{The Bony's paraproduct}

Next, we follow the Bony's paraproduct decomposition. Given tempered distributions $u,v\in S'(\mathbb{R}^{d})$, it is valid the decomposition
\begin{equation*}
 vw=T_{v}w+T_{w}v+R(v,w),   
\end{equation*}
where the Bony's paraproduct operator $T_v(\cdot)$ is given by 
\begin{equation*}
T_{v}(w)=\displaystyle{\sum_{l\in\mathbb{Z}}}\dot{S}_{l-1}v \dot{\Delta}_{l}w,   
\end{equation*}
and the remainder part $R(\cdot,\cdot)$ is given by 
\begin{equation*}
R(v,w)=\displaystyle{\sum_{l\in\mathbb{Z}}}\dot{\Delta}_{l}v\tilde{\dot{\Delta}}_{l}w,\,\,\mbox{where},\,\,
\tilde{\dot{\Delta}}_{l}w=\displaystyle{\sum_{\mid l-l'\mid\leq 1}}\dot{\Delta}_{l'}w.
\end{equation*}
Moreover, we also denote $\tilde{\phi}_{l}=\phi_{l-1}+\phi_{l}+\phi_{l+1}$. For further information on 
Bony's paraproduct see (\cite{Lemarie-NS-century-2018}, Chapter 3), (\cite{Bahouri-Chemin-Danchin-Fourier-PDE-2011}, 
Section 2.6) and the original work by J. M. Bony \cite{Bony-symbolique-1981}.

\subsection{An abstract lemma}

In this section we recall an abstract lemma related to a fixed point argument (\cite{Lemarie-NS-century-2018}, Theorem 13.2).

First, we note that the mild formulation has the abstract form 
\begin{equation}\label{mild-formulation-abstract}
u(t)=z_{0}+B(u,u)(t),     
\end{equation}
for all $t\in [0,\infty)$, where $z_{0}=S_{\nu}(t)u_{0}$ and $B(\cdot,\cdot)$ is the given bilinear operator. 

The following lemma is the fixed point scheme which we shall apply in order to solve the abstract equation associated to the mild formulation for the Navier-Stokes system.
\begin{lema}\label{fixed-point-scheme}
Let $\left(Z,\parallel \cdot\parallel\right)$ be a Banach space and $B:Z\times Z\rightarrow Z$ be a bilinear operator for which there is $K>0$ such that $\parallel B\parallel_{{\cal B}(Z)}\leq K$. Then, for $0<\varepsilon<\frac{1}{4K}$ and $\parallel z_{0}\parallel_{Z}\leq \varepsilon$, the abstract equation $z=z_{0}+B(z,z)$ has a unique solution in the closed ball $B_{2\varepsilon}(0)=\{z\in Z;\parallel z\parallel_{Z}\leq 2\varepsilon\}$ and $\parallel z\parallel_{Z}\leq 2\parallel z_{0}\parallel_{Z}$. Moreover, the solution depends continuously on initial data, in the following sense: for $\tilde{z}_{0}\in Z$, with $\parallel\tilde{z}_{0}\parallel_{Z}\leq \varepsilon$, if $\tilde{z}=\tilde{z}_{0}+B(\tilde{z},\tilde{z})$ and $\parallel \tilde{z}\parallel_{Z}\leq 2\varepsilon$, then 
\begin{equation*}
 \parallel z-\tilde{z}\parallel_{Z}\leq \left(1-4K\varepsilon\right)^{-1}\parallel z_{0}-\tilde{z}_{0}\parallel_{Z}.   
\end{equation*}
\end{lema}
\subsection{The linear estimates}
In this section we prove the main linear estimates related to mixed-norm Besov-Lebesgue spaces and mixed-norm Fourier-Besov-Lebesgue spaces, involving the Stokes semigroup.
\begin{lema}\label{estimates-lin-besov}
 Let us take $\vec{p}\in[1,\infty]^{d}$ and $q\in[1,\infty]$. For any $\sigma\in\mathbb{R}$, there exists a constant $C>0$ independent on $\nu$, such that
\begin{equation*}
\begin{split}
 (i.) &\parallel S_{\nu}(\cdot)u_{0}\parallel_{{\cal L}^{\infty}\left(I;\dot{{\cal B}}^{\sigma}_{\vec{p},q}(\mathbb{R}^{d})\right)}\leq C \parallel u_{0}\parallel_{\dot{{\cal B}}^{\sigma}_{\vec{p},q}(\mathbb{R}^{d})}, 
 \\
 (ii.)&\parallel S_{\nu}(\cdot)u_{0}\parallel_{{\cal L}^{1}\left(I;\dot{{\cal B}}^{\sigma+2}_{\vec{p},q}(\mathbb{R}^{d})\right)}\leq \frac{C}{\nu} \parallel u_{0}\parallel_{\dot{{\cal B}}^{\sigma}_{\vec{p},q}(\mathbb{R}^{d})},  
\end{split}
\end{equation*}
 for any $u_{0}\in \dot{{\cal B}}^{\sigma}_{\vec{p},q}(\mathbb{R}^{d})$.
\end{lema}
\begin{proof}
\begin{itemize}
\item [1.] Proof of item (i.): From $e^{-ct2^{2l}}\leq 1$, for all $t>0$ and any $l\in\mathbb{Z}$ we obtain (i.) if we consider the definition of the norm on ${\cal L}^{\infty}\left(I;{\cal B}^{\sigma}_{\vec{p},q}(\mathbb{R}^{d})\right)$. 
\item [2.] Proof of item (ii.): Since $\parallel e^{\nu t\Delta}\dot{\Delta}_{l}u\parallel_{L^{\vec{p}}(\mathbb{R}^{d})}\leq e^{-c\nu t2^{2l}} \parallel\dot{\Delta}_{l}u\parallel_{L^{\vec{p}}(\mathbb{R}^{d})}$, and the definition of the norm on  ${\cal L}^{1}\left(I;{\cal B}^{\sigma+2}_{\vec{p},q}(\mathbb{R}^{d})\right)$, we get
\begin{equation*}
 \begin{split}
 \parallel S_{\nu}u_{0}\parallel_{{\cal L}^{1}\left(I;{\cal B}^{\sigma+2}_{\vec{p},q}\right)}&\leq
 \parallel 2^{l(\sigma+2)}\displaystyle{\int_{0}^{\infty}}\parallel e^{\nu t'\Delta}\dot{\Delta}_{l}u_{0}\parallel_{L^{\vec{p}}(\mathbb{R}^{d})}\,dt'\parallel_{l^{q}(\mathbb{Z})}\\
 &\leq
 \parallel 2^{l(\sigma+2)}\left(\displaystyle{\int_{0}^{\infty}}e^{-c\nu t2^{2l}}\,dt'\right)\parallel\dot{\Delta}_{l}u_{0}\parallel_{L^{\vec{p}}(\mathbb{R}^{d})}\parallel_{l^{q}(\mathbb{Z})},
 \end{split}   
\end{equation*}
which implies (ii.) since the integral is bounded by $\frac{C}{\nu}2^{-2l}$.
\end{itemize}
\end{proof}

\begin{lema}
Let us consider $\vec{p}\in[1,\infty]^{d}$ and $q\in[1,\infty]$. For each $s\in \mathbb{R}$, there is a constant $C>0$ independent of $\nu$, such that
\begin{equation*}
\begin{split}
(i.)&\parallel S_{\nu}(\cdot)u_{0}\parallel_{{\cal L}^{\infty}\left(I;\dot{{\cal FB}}^{s}_{\vec{p},q}(\mathbb{R}^{d})\right)}\leq C\parallel u_{0}\parallel_{\dot{{\cal FB}}^{s}_{\vec{p},q}(\mathbb{R}^{d})},\\    
(ii.)&
\parallel S_{\nu}(\cdot)u_{0}\parallel_{{\cal L}^{1}\left(I;\dot{{\cal FB}}^{s+2}_{\vec{p},q}(\mathbb{R}^{d})\right)}\leq \frac{C}{\nu}\parallel u_{0}\parallel_{\dot{{\cal FB}}^{s}_{\vec{p},q}(\mathbb{R}^{d})},    
\end{split}
\end{equation*} 
for every $u_{0}\in \dot{{\cal FB}}^{s}_{\vec{p},q}(\mathbb{R}^{d})$.
\end{lema}

\begin{proof}
\begin{itemize}
\item [1.]Proof of item (i.): Since $e^{-\nu t\mid\eta\mid^{2}}\leq 1$, for all $t\geq 0$ and $\eta\in\mathbb{R}^{d}$, we easily get $(i)$, by considering the definition of ${\cal L}^{\infty}\left(I;\dot{{\cal FB}}^{s}_{\vec{p},q}(\mathbb{R}^{d})\right)$.  
\item [2.] Proof of item (ii.): By definition of the norm on ${\cal L}^{1}\left(I;\dot{{\cal FB}}^{s+2}_{\vec{p},q}(\mathbb{R}^{d})\right)$, we get
\begin{equation*}
 \begin{split}
  \parallel S_{\nu}(\cdot)u_{0}\parallel_{{\cal L}^{1}\left(I;\dot{{\cal FB}}^{s+2}_{\vec{p},q}(\mathbb{R}^{d})\right)}&\leq 
  \displaystyle{\parallel 
  2^{l(s+2)}\int_{0}^{\infty}\parallel e^{-\nu t'\mid \eta\mid^{2}}\phi_{l}(\eta)\hat{u}_{0}\parallel_{L^{\vec{p}}(\mathbb{R}^{d})}\,dt'
  \parallel}_{l^{q}(\mathbb{Z})}\\
  &\leq 
  \displaystyle{\parallel 
  2^{l(s+2)}\left(\int_{0}^{\infty} e^{-\nu t'2^{2l}}\,dt'\right)\parallel \phi_{l}(\eta)\hat{u}_{0}\parallel_{L^{\vec{p}}(\mathbb{R}^{d})}
  \parallel}_{l^{q}(\mathbb{Z})}.
 \end{split}   
\end{equation*}
Since the last integral is bounded by $\frac{C
}{\nu}2^{-2l}$, we have the desired $(ii.)$.
\end{itemize}
\end{proof}

\section{The bilinear estimates}
\label{bilinear-estimates-s}

In this section, we first introduce an auxiliary linear operator and we establish 
some basic estimates for this operator. We also set the precise estimate for 
the bilinear term associated to the mild formulation. 

\medskip
First, let us define the following auxiliary linear operator
\begin{equation}\label{auxiliar-linearop}
 A(g)(t,\cdot):= \displaystyle{\int_{0}^{t}}e^{\nu(t-t')\Delta}\mathbb{P}g(t',\cdot)\,dt'   
\end{equation}
which in Fourier variables is given by
\begin{equation*}
 [A(g)]^{\wedge}(\eta,t)=\displaystyle{\int_{0}^{t}e^{-\nu(t-t')\mid\eta\mid^{2}}\hat{\mathbb{P}}(\eta)\hat{g}(\eta,t')\,dt'}.   
\end{equation*}
Then, we have the following 
\begin{lema}
Let us take $\vec{p}\in[1,\infty]^{d}$, $q\in[1,\infty]$ and $s\in\mathbb{R}$. Then, there exists $C'>0$ independent on $\nu$ such that
\begin{equation*}
\begin{split}
    (i.)&\parallel A(g)\parallel_{{\cal L}^{\infty}\left(I;\dot{{\cal B}}^{\sigma}_{\vec{p},q}(\mathbb{R}^{d})\right)}\leq C'
    \parallel g\parallel_{{\cal L}^{1}\left(I;\dot{{\cal B}}^{\sigma}_{\vec{p},q}(\mathbb{R}^{d})\right)},
    \\
(ii.)& \parallel A(g)\parallel_{{\cal L}^{1}\left(I;\dot{{\cal B}}^{\sigma+2}_{\vec{p},q}(\mathbb{R}^{d})\right)}\leq \frac{C'}{\nu}
    \parallel g\parallel_{{\cal L}^{1}\left(I;\dot{{\cal B}}^{\sigma}_{\vec{p},q}(\mathbb{R}^{d})\right)},
\end{split}
\end{equation*}
for all $g\in {\cal L}^{1}\left(I;\dot{{\cal B}}^{\sigma}_{\vec{p},q}(\mathbb{R}^{d})\right)$.
\end{lema}
\begin{proof}
\begin{itemize}
First, since the Riesz transform ${\cal R}_{j}$ is given by ${\cal R}_{j}[\dot{\Delta}_{l}g]={\cal F}^{-1}[\frac{i\eta_{j}}{\mid \eta\mid} (\dot{\Delta}_{l}g)^{\wedge}(\eta)]$, due to the Young's inequality we conclude that 
    \begin{equation*}
     \parallel {\cal R}_{j}\dot{{\Delta}}_{l}g\parallel_{L^{\vec{p}}(\mathbb{R}^{d})}\leq C \parallel \dot{{\Delta}}_{l}g\parallel_{L^{\vec{p}}(\mathbb{R}^{d})},  
    \end{equation*}
    for some constant $C>0$ independent on $j, l$ and $g$. Therefore, by definition of the Leray operator $\mathbb{P}$, we obtain
    \begin{equation*}
     \parallel \mathbb{P}\dot{{\Delta}}_{l}g\parallel_{L^{\vec{p}}(\mathbb{R}^{d})}\leq C \parallel\dot{{\Delta}}_{l}g\parallel_{L^{\vec{p}}(\mathbb{R}^{d})},  
    \end{equation*}
    for some constant $C>0$ independent on $l$ and $g$.
    \item [1.] Proof of item (i.): Since $\parallel e^{\nu t\Delta}\dot{\Delta}_{l}u\parallel_{L^{\vec{p}}(\mathbb{R}^{d})}\leq e^{-\nu t2^{2l}}\parallel \dot{\Delta}_{l}u\parallel_{L^{\vec{p}}(\mathbb{R}^{d})}$ and from the definition of the norm on  ${\cal L}^{\infty}\left(I;\dot{{\cal B}}^{\sigma}_{\vec{p},q}(\mathbb{R}^{d})\right)$, we get
\begin{equation*}
\begin{split}
   \parallel A(g)\parallel_{{\cal L}^{\infty}\left(I;\dot{{\cal B}}^{\sigma}_{\vec{p},q}(\mathbb{R}^{d})\right)}&\leq \parallel 2^{l\sigma}\parallel \displaystyle{\int_{0}^{t}}e^{-\nu c(t-t')2^{2l}}\parallel \dot{\Delta}_{l}g\parallel_{L^{\vec{p}}(\mathbb{R}^{d})}\,dt'\parallel_{L^{\infty}(I)}\parallel_{l^{q}(\mathbb{Z})}\\
   &\leq \parallel 2^{l\sigma}\parallel \parallel\dot{\Delta}_{l}g\parallel_{L^{1}\left((0,t);L^{\vec{p}}(\mathbb{R}^{d})\right)}\parallel_{L^{\infty}(I)}\parallel_{l^{q}(\mathbb{Z})}
\end{split}    
\end{equation*}
and this implies (i.) since $(0,t)\subset I$.
\item [2.] Proof of item (ii.): Again, by the definition of the norm on ${\cal L}^{1}\left(I;\dot{{\cal B}}^{\sigma+2}_{\vec{p},q}(\mathbb{R}^{d})\right)$, we have 
\begin{equation*}
 \begin{split}
  \parallel A(g)\parallel_{{\cal L}^{1}\left(I;\dot{{\cal B}}^{\sigma+2}_{\vec{p},q}(\mathbb{R}^{d})\right)}&\leq 
  \parallel 2^{l(\sigma+2)}\parallel \displaystyle{\int_{0}^{t}}e^{-c\nu(t-t')2^{2l}}\parallel \dot{\Delta}_{l}g\parallel_{L^{\vec{p}}(\mathbb{R}^{d})}\,dt'\parallel_{L^{1}(I)}\parallel_{l^{q}(\mathbb{Z})}\\
  &\leq \parallel 2^{l(\sigma+2)}\displaystyle{\int_{0}^{\infty}}
  \left(\displaystyle{\int_{t'}^{\infty}}e^{-c\nu(t-t')2^{2l}}\,dt\right)\cdot \parallel \dot{\Delta}_{l}g(t')\parallel_{L^{\vec{p}}(\mathbb{R}^{d})}\,dt'
 \end{split}   
\end{equation*}
and this implies (ii.) since the integral is bounded by $\frac{C}{\nu}2^{-2l}$.
\end{itemize}
\end{proof}
\begin{obs}
One observes that, 
the Riesz transform is bounded on homogeneous Besov mixed-Lebesgue spaces. 
In fact, the proof is similar to the usual homogeneous Besov spaces (\cite{Sawano-homogeneous-2020}, Theorem 1.3). 
This assertion implies, by definition, that the Leray projector $\mathbb{P}$ is also bounded on these homogeneous Besov mixed-Lebesgue spaces.   
\end{obs}
\begin{lema}
Let us consider $\vec{p}\in [1,\infty]^{d}$, $q\in[1,\infty]$ and $s\in\mathbb{R}$. Then, there exists $C'>0$ independent on $\nu$ such that
\begin{equation*}
\begin{split}
(i.)&\parallel A(g)\parallel_{{\cal L}^{\infty}\left(I;\dot{{\cal FB}}^{s}_{\vec{p},q}(\mathbb{R}^{d})\right)}\leq C'  
\parallel g\parallel_{{\cal L}^{1}(I;\dot{{\cal FB}}^{s}_{\vec{p},q}(\mathbb{R}^{d}))},\\
(ii.)&\parallel A(g)\parallel_{{\cal L}^{1}\left(I;\dot{{\cal FB}}^{s+2}_{\vec{p},q}\right)}\leq \frac{C'}{\nu}  \cdot
\parallel g\parallel_{{\cal L}^{1}\left(I;\dot{{\cal FB}}^{s}_{\vec{p},q}(\mathbb{R}^{d})\right)},    
\end{split}
\end{equation*}
for all $g\in {\cal L}^{1}\left(I;\dot{{\cal FB}}^{s}_{\vec{p},q}(\mathbb{R}^{d})\right)$.
\end{lema}
\begin{proof}
First, recall that $\supp \phi_{l}\subset 2^{l}{\cal A}=\{\eta\in\mathbb{R}^{d};\frac{3}{4}\cdot 2^{l}\leq \mid \eta\mid\leq \frac{8}{3}\cdot 2^{l}\}$ and ${\cal R}_{j}[\Delta_{l}g]={\cal F}^{-1}[\frac{i\eta_{j}}{\mid \eta\mid}\phi_{l}\hat{g}]$. Then $\mid [{\cal R}_{j}(\Delta_{l}g)]^{\wedge}(\eta)\mid=\mid \frac{i\eta_{j}}{\mid\eta\mid}\mid \mid\phi_{l}\hat{g}\mid\leq \mid\phi_{l}\hat{g}\mid$, for all $\eta\in\mathbb{R}^{d}-\{0\}$. This implies that $\mid {\cal F}[\mathbb{P}\Delta_{l}g](\eta)\mid\leq \mid \phi_{l}\hat{g}(\eta)\mid$, for all $\eta\in\mathbb{R}^{d}-\{0\}$.
\begin{itemize}
\item [(i.)] We note that, by definition of the norm on ${\cal L}^{\infty}\left(I;\dot{{\cal FB}}^{s}_{\vec{p},q}(\mathbb{R}^{d})\right)$, we have
\begin{equation*}
\begin{split}
\parallel A(g)\parallel_{{\cal L}^{\infty}\left(I;\dot{{\cal FB}}^{s}_{\vec{p},q}(\mathbb{R}^{d})\right)}&\leq
\parallel 2^{ls}\parallel \displaystyle{\int_{0}^{t}}
e^{-\nu(t-t')2^{2l}}\parallel \phi_{l}\hat{g}\parallel_{L^{\vec{p}}(\mathbb{R}^{d})}\,dt'\parallel_{L^{\infty}(I)}\parallel_{l^{q}(\mathbb{Z})}\\
&\leq \parallel 2^{ls}\parallel \parallel \phi_{l}\hat{g}\parallel_{L^{1}\left((0,t);L^{\vec{p}}(\mathbb{R}^{d})\right)}\parallel_{L^{\infty}(I)}\parallel_{l^{q}(\mathbb{Z})}
\end{split}    
\end{equation*}
which implies iten (i.) since $(0,t)\subset I$.
\item [(ii.)] Again, by definition of the norm on ${\cal L}^{1}\left(I;\dot{{\cal FB}}^{s+2}_{\vec{p},q}(\mathbb{R}^{d})\right)$, we have
\begin{equation*}
\begin{split}
\parallel A(g)\parallel_{{\cal L}^{1}\left(I;\dot{{\cal FB}}^{s+2}_{\vec{p},q}(\mathbb{R}^{d})\right)}&\leq
\parallel 2^{l(s+2)}\parallel \displaystyle{\int_{0}^{t}}
e^{-\nu(t-t')2^{2l}}\parallel \phi_{l}\hat{g}\parallel_{L^{\vec{p}}(\mathbb{R}^{d})}\,dt'\parallel_{L^{1}(I)}\parallel_{l^{q}(\mathbb{Z})}\\
&\leq \parallel 2^{l(s+2)}\cdot\frac{C'}{\nu}\cdot 2^{-2l} \parallel \phi_{l}\hat{g}\parallel_{L^{1}\left(I;L^{\vec{p}}(\mathbb{R}^{d})\right)}\parallel_{l^{q}(\mathbb{Z})}
\end{split}    
\end{equation*}
which implies iten (ii.), by considering the sum of the exponents of $2$.
\end{itemize}
\end{proof}
\begin{lema}\label{estimates-bilinear-besov}
 Let us consider $\vec{p}\in[1,\infty]^{d}$ such that $\displaystyle{\sum_{i=1}^{d}}\frac{1}{p_{i}}>0$, $q\in[1,\infty]$ and $\sigma=-1+\displaystyle{\sum_{i=1}^{d}}\frac{1}{p_{i}}$. Then, there exists a constant $K_{0}>0$ independent on $\nu$ such that
 \begin{equation*}
  \parallel B(v,w)\parallel_{Z}\leq K_{0}\max{(1,\frac{1}{\nu})}\parallel v\parallel_{Z}\parallel w\parallel_{Z},   
 \end{equation*}
 for each $v,w\in Z={\cal L}^{\infty}\left(I;\dot{{\cal B}}^{\sigma}_{\vec{p},q}(\mathbb{R}^{d})\right)\bigcap {\cal L}^{1}\left(I;\dot{{\cal B}}^{\sigma+2}_{\vec{p},q}(\mathbb{R}^{d})\right)$.
\end{lema}
\begin{proof}
First, we observe that since $B(v,w)(t,\cdot)=A(\mbox{div}(v\otimes w))(t,\cdot)$, where the linear auxiliar operator $A(\cdot)$ is given in (\ref{auxiliar-linearop}), we get 
$$
\begin{aligned}
 \parallel B(v,w)\parallel_{Z}&\leq c\max{(1,\frac{1}{\nu})} \parallel \mbox{div}\left(v\otimes w\right)\parallel_{{\cal L}^{1}(I;\dot{{\cal B}}^{\sigma}_{\vec{p},q}(\mathbb{R}^{d}))}
 \\
& \leq c\max{(1,\frac{1}{\nu})}\parallel v\otimes w\parallel_{{\cal L}^{1}(I;\dot{{\cal B}}^{\sigma +1}_{\vec{p},q}(\mathbb{R}^{d}))},   
\end{aligned}
$$
where $c$ is a constant which does not depend on $\nu$. Therefore, it is enought o prove that 
\begin{equation*}
 \parallel v\otimes w\parallel_{{\cal L}^{1}(I;\dot{{\cal B}}^{\sigma +1}_{\vec{p},q}(\mathbb{R}^{d}))}\leq \parallel v\parallel_{Z}\parallel w\parallel_{Z}.
\end{equation*}
In order to prove this estimate we consider the Bony's paraproduct decomposition
\begin{equation*}
 \begin{split}
  \dot{\Delta}_{l}\left(vw\right)&=\displaystyle{\sum_{\mid l-l'\mid\leq 4}}\dot{\Delta}_{l}\left( \dot{S}_{l'-1}v\dot{\Delta}_{l'}w\right) + \displaystyle{\sum_{\mid l-l'\mid\leq 4}}\dot{\Delta}_{l}\left( \dot{S}_{l'-1}w\dot{\Delta}_{l'}v\right) + \displaystyle{\sum_{l'\geq l-2}}\dot{\Delta}_{l}\left(\dot{\Delta}_{l'}v\tilde{\dot{\Delta}}_{l'}w\right)\\
  &=I^{l}_{1} + I^{l}_{2} + I^{l}_{3}.
 \end{split}   
\end{equation*}
Estimates for $I^{l}_{1}$: We first observe that the Bernstein's type inequality Lemma 2.4 (i.) (\ref{properties-Besov-spaces-mixed}) implies $\parallel \dot{\Delta}_{l''}v\parallel_{L^{\infty}(\mathbb{R}^{d})}\leq 2^{l''(\sigma +1)}\parallel \dot{\Delta}_{l''}v\parallel_{L^{\vec{p}}(\mathbb{R}^{d})}$ and thus we have
\begin{equation*}
\begin{split}
\parallel I^{l}_{1}\parallel_{L^{1}\left(I;L^{\vec{p}}(\mathbb{R}^{d})\right)}&\leq 
\displaystyle{\sum_{\mid l-l'\mid\leq 4}}\parallel \parallel \dot{S}_{l'-1}v\dot{\Delta}_{l'}w\parallel_{L^{\vec{p}}(\mathbb{R}^{d})}\parallel_{L^{1}(I)}\\
&\leq \displaystyle{\sum_{\mid l-l'\mid\leq 4}}
\parallel \displaystyle{\sum_{l''\leq l'}}\parallel\dot{\Delta}_{l''}v\parallel_{L^{\infty}(\mathbb{R}^{d})}\parallel\dot{\Delta}_{l'}w\parallel_{L^{\vec{p}}(\mathbb{R}^{d})}\parallel_{L^{1}(I)}\\
&\leq \displaystyle{\sum_{\mid l-l'\mid\leq 4}} \left(\displaystyle{\sum_{l''\leq l'}}
2^{\sigma +1}\parallel \dot{\Delta}_{l''}v\parallel_{L^{\infty}\left(I;L^{\vec{p}}(\mathbb{R}^{d})\right)}\right) \parallel \dot{\Delta}_{l'}w\parallel_{L^{1}\left(I;L^{\vec{p}}(\mathbb{R}^{d})\right)}\\
&= \displaystyle{\sum_{\mid l-l'\mid\leq 4}} \left(\displaystyle{\sum_{l''\leq l'}}
2^{l''}\cdot 2^{l''\sigma}\parallel \dot{\Delta}_{l''}v\parallel_{L^{\infty}\left(I; L^{\vec{p}}(\mathbb{R}^{d})\right)}\right) \parallel \dot{\Delta}_{l'}w\parallel_{L^{1}\left(I;L^{\vec{p}}(\mathbb{R}^{d})\right)}.
\end{split}   
\end{equation*}
By H\"older's inequality for series, we get
\begin{equation*}
 \begin{split}
 \parallel I^{l}_{1}\parallel_{L^{1}\left(I;L^{\vec{p}}(\mathbb{R}^{d})\right)}&\leq    
  \displaystyle{\sum_{\mid l-l'\mid\leq 4}}\left(\displaystyle{\sum_{l''\leq l'}2^{l''q'}}\right)^{\frac{1}{q'}}\cdot \parallel \dot{\Delta}_{l'}w\parallel_{L^{1}\left(I;L^{\vec{p}}(\mathbb{R}^{d})\right)}\parallel v \parallel_{{\cal L}^{\infty}\left(I;\dot{{\cal B}}^{\sigma}_{\vec{p},q}(\mathbb{R}^{d})\right)}\\
  &\leq \displaystyle{\sum_{\mid l-l'\mid\leq 4}}2^{l'} \parallel \dot{\Delta}_{l'}w\parallel_{L^{1}\left(I;L^{\vec{p}}(\mathbb{R}^{d})\right)}\parallel v\parallel_{Z}\\
&=  \displaystyle{\sum_{\mid l-l'\mid\leq 4}}2^{-l'(\sigma+1)}\cdot 2^{l'(\sigma+2)}\parallel \dot{\Delta}_{l'}w\parallel_{L^{1}\left(I;L^{\vec{p}}(\mathbb{R}^{d})\right)}\parallel v\parallel_{Z}.
 \end{split}   
\end{equation*}
Multiplying by $2^{l(\sigma+1)}$ in both sides, we have
\begin{equation*}
 \begin{split}
  2^{l(\sigma+1)}\parallel I^{l}_{1}\parallel_{L^{1}\left(I;L^{\vec{p}}(\mathbb{R}^{d})\right)}&\leq
  \displaystyle{\sum_{l'}}2^{(l-l')(\sigma+1)}\chi_{\{k;\mid k\mid\leq 4\}}\cdot \left(2^{l'(\sigma+2)}\parallel \dot{\Delta}_{l'}w\parallel_{L^{1}\left(I;L^{\vec{p}}(\mathbb{R}^{d})\right)}\right)\cdot\parallel v\parallel_{Z}\\
  &=\left(a_{k}\ast b_{l'}\right) _{l}\cdot\parallel v\parallel_{Z},
 \end{split}
\end{equation*}
where $a_{k}=2^{k(\sigma+1)}\chi_{\{k;\mid k\mid\leq 4\}}$ and $b_{l'}=2^{l'(\sigma+2)}\parallel\dot{\Delta}_{l'}w\parallel_{L^{1}\left(I;L^{\vec{p}}(\mathbb{R}^{d})\right)}$. By applying the Young's inequality for series, we have
\begin{equation*}
\begin{split}
 \parallel  2^{l(\sigma+1)}\parallel I^{l}_{1}\parallel_{L^{1}\left(I;L^{\vec{p}}(\mathbb{R}^{d})\right)}\parallel_{l^{q}(\mathbb{Z})}&\leq
 \parallel a_{k}\parallel_{l^{1}(\mathbb{Z})}\parallel b_{l'}\parallel_{l^{q}(\mathbb{Z})}\parallel v\parallel_{Z}\\
 &\leq \parallel w\parallel_{{\cal L}^{1}\left(I;\dot{{\cal B}}^{\sigma+2}_{\vec{p},q}(\mathbb{R}^{d})\right)} \parallel v\parallel_{Z},
\end{split}    
\end{equation*}
which implies the desired estimate for $I^{l}_{1}$. 
A similar procedure is enough to prove that
\begin{equation*}
\parallel  2^{l(\sigma+1)}\parallel I^{l}_{2}\parallel_{L^{1}\left(I;L^{\vec{p}}(\mathbb{R}^{d})\right)}\parallel_{l^{q}(\mathbb{Z})}\leq \parallel w\parallel_{Z}\parallel v\parallel_{Z}.    
\end{equation*}
Estimates for the term $I^{l}_{3}$: By the Bernstein's type inequality Lemma 2.4 (i.) (\ref{properties-Besov-spaces-mixed}), we obtain
\begin{equation*}
 \begin{split}
 \parallel I^{l}_{3}\parallel_{L^{1}\left(I;L^{\vec{p}}(\mathbb{R}^{d})\right)}&\leq 
 \displaystyle{\sum_{l'\geq l-2}}\parallel \dot{\Delta}_{l'}v\parallel_{L^{1}\left(I;L^{\vec{p}}(\mathbb{R}^{d})\right)} \left(
 \displaystyle{\sum_{\mid l''-l'\mid\leq 1}}2^{l''(\sigma+1))}\parallel \dot{\Delta}_{l''}w\parallel_{L^{\infty}\left(I;L^{\vec{p}}(\mathbb{R}^{d})\right)}\right)\\
 &=\displaystyle{\sum_{l'\geq l-2}}\parallel \dot{\Delta}_{l'}v\parallel_{L^{1}\left(I;L^{\vec{p}}(\mathbb{R}^{d})\right)} \left(
 \displaystyle{\sum_{\mid l''-l'\mid\leq 1}}2^{l''}\cdot 2^{l''\sigma}\parallel \dot{\Delta}_{l''}w\parallel_{L^{\infty}\left(I;L^{\vec{p}}(\mathbb{R}^{d})\right)}
 \right).
 \end{split}   
\end{equation*}
By the H\"older's inequality for series we obtain
\begin{equation*}
 \begin{split}
 \parallel I^{l}_{3}\parallel_{L^{1}\left(I;L^{\vec{p}}(\mathbb{R}^{d})\right)}&\leq 
 \displaystyle{\sum_{l'\geq l-2}}\parallel \dot{\Delta}_{l'}v\parallel_{L^{1}\left(I;L^{\vec{p}}(\mathbb{R}^{d})\right)} \left(\displaystyle{\sum_{\mid l''-l'\mid \leq 1}} 2^{l''q'}\right)^{\frac{1}{q'}}\parallel w\parallel_{{\cal L}^{\infty}\left(I;\dot{{\cal B}}^{\sigma}_{\vec{p},q}\right)}\\
&=\displaystyle{\sum_{l'\geq l-2}}\parallel \dot{\Delta}_{l'}v\parallel_{L^{1}\left(I;L^{\vec{p}}(\mathbb{R}^{d})\right)} 2^{l'} \parallel w\parallel_{{\cal L}^{\infty}\left(I;\dot{{\cal B}}^{\sigma}_{\vec{p},q}\right)}\\
&=\displaystyle{\sum_{l'\geq l-2}}2^{-l'(\sigma+1)}\cdot 2^{l'(\sigma+2)}\parallel \dot{\Delta}_{l'}v\parallel_{L^{1}\left(I;L^{\vec{p}}(\mathbb{R}^{d})\right)}\parallel w\parallel_{Z}.
 \end{split}   
\end{equation*}
Multiplying by $2^{l(\sigma+1)}$ we obtain
\begin{equation*}
 \begin{split}
  2^{l(\sigma+1)} \parallel I^{l}_{3}\parallel_{L^{1}\left(I;L^{\vec{p}}(\mathbb{R}^{d})\right)}&\leq \displaystyle{\sum_{l'}}
  2^{(l-l')(\sigma+1)}\chi_{\{k; k\leq 2\}}(l-l')\cdot 2^{l'(\sigma+2)}\parallel \dot{\Delta}_{l'}v\parallel_{L^{1}\left(I;L^{\vec{p}}(\mathbb{R}^{d})\right)}\parallel w\parallel_{Z}\\
  &=\left(a_{k}\ast b_{l'}\right)_{l},
 \end{split}   
\end{equation*}
where $a_{k}=2^{k(\sigma+1)}\chi_{\{k;k\leq 2\}}$ and $b_{l'}=2^{l'(\sigma+2)}\parallel \dot{\Delta}_{l'}v\parallel_{L^{1}\left(I;L^{\vec{p}}(\mathbb{R}^{d})\right)}$.
Now we apply the Young's inequality for series to get
\begin{equation*}
\begin{split}
 \parallel   2^{l(\sigma+1)} \parallel I^{l}_{3}\parallel_{L^{1}\left(I;L^{\vec{p}}(\mathbb{R}^{d})\right)}\parallel_{l^{q}(\mathbb{Z})}&\leq 
 \parallel a_{k}\parallel_{l^{1}(\mathbb{Z})}\parallel b_{l'}\parallel_{l^{q}(\mathbb{Z})}\parallel w\parallel_{Z}\\
 &\leq \parallel v\parallel_{{\cal L}^{1}\left(I;\dot{{\cal B}}^{\sigma+2}_{\vec{p},q}\right)}\parallel w\parallel_{Z},
 \end{split}
 \end{equation*}
which implies the desired estimate, since the sequence $(a_{k})\in l^{1}(\mathbb{Z})$ whenever $\sigma+1>0$.
 
\end{proof}
\begin{lema}
 Let us consider $\vec{p}\in[1,\infty]^{d}$, such that $\displaystyle{\sum_{i=1}^{d}}\left(1-\frac{1}{p_{i}}\right)>0$, $q\in[1,\infty]$ and $s=-1+\displaystyle{\sum_{i=1}^{d}}\left(1-\frac{1}{p_{i}}\right)$. Then, there exists a constant $K_{0}>0$ independent on $\nu$ such that
 \begin{equation*}
  \parallel B(v,w)\parallel_{Z}\leq K_{0}\max{(1,\frac{1}{\nu})}\parallel v\parallel_{Z}\parallel w\parallel_{Z},   
 \end{equation*}
 for each $v,w\in Z={\cal L}^{\infty}\left(I;\dot{{\cal FB}}^{s}_{\vec{p},q}(\mathbb{R}^{d})\right)\bigcap {\cal L}^{1}\left(I;\dot{{\cal FB}}^{s+2}_{\vec{p},q}(\mathbb{R}^{d})\right)$.
\end{lema}
\begin{proof}
First, we observe that since $B(v,w)(t,\cdot)=A(\mbox{div}(v\otimes w))(t,\cdot)$, where the linear auxiliar operator $A(\cdot)$ is given in (\ref{auxiliar-linearop}), we get 
\begin{equation*}
 \parallel B(v,w)\parallel_{Z}\leq c\max{(1,\frac{1}{\nu})}\parallel \mbox{div}\left(v\otimes w\right)\parallel_{{\cal L}^{1}(I;\dot{{\cal FB}}^{s}_{\vec{p},q}(\mathbb{R}^{d}))}\leq c\max{(1,\frac{1}{\nu})}\parallel v\otimes w\parallel_{{\cal L}^{1}(I;\dot{{\cal FB}}^{s +1}_{\vec{p},q}(\mathbb{R}^{d}))},   
\end{equation*}
where $c$ is a constant which does not depend on $\nu$. Therefore, it is enought o prove that 
\begin{equation*}
 \parallel v\otimes w\parallel_{{\cal L}^{1}(I;\dot{{\cal FB}}^{s +1}_{\vec{p},q}(\mathbb{R}^{d}))}\leq \parallel v\parallel_{Z}\parallel w\parallel_{Z}.
\end{equation*}
By the Bony's decomposition we have
\begin{equation*}
\begin{split}
\phi_{l}\left(vw\right)^{\wedge}&=\displaystyle{\sum_{\mid l-l'\mid\leq 4}}\phi_{l}[\left(\dot{S}_{l'-1}v\right)^{\wedge}\ast\left(\phi_{l'}\hat{w}\right)] + \displaystyle{\sum_{\mid l-l'\mid\leq 4}}\phi_{l}[\left(\dot{S}_{l'-1}w\right)^{\wedge}\ast\left(\phi_{l'}\hat{v}\right)] + \displaystyle{\sum_{l'\geq l-2}}\phi_{l}[\left(\phi_{l'}\hat{v}\right)\ast\left( \tilde{\phi}_{l'}\hat{w}\right)]\\
&=I^{l}_{1}+I^{l}_{2}+I^{l}_{3}.
\end{split}
\end{equation*}
Estimates for $I^{l}_{1}$: By the Bernstein's type inequality Lemma 2.6 (i.) (\ref{properties-Fourier-Besov-mixed}), we have $\parallel\phi_{l''}\hat{v}\parallel_{L^{1}(\mathbb{R}^{d})}\leq 2^{l''(s+1)}\parallel\phi_{l''}\hat{v}\parallel_{L^{\vec{p}}(\mathbb{R}^{d})}$ we have
\begin{equation*}.
\begin{split}
\parallel I^{l}_{1}\parallel_{L^{1}\left(I;L^{\vec{p}}(\mathbb{R}^{d})\right)}&\leq \displaystyle{\sum_{\mid l-l'\mid\leq 4}}\parallel \parallel \left(S_{l'-1}v\right)^{\wedge}\ast \left(\phi_{l'}\hat{w}\right)\parallel_{L^{\vec{p}}(\mathbb{R}^{d})}\parallel_{L^{1}(I)}\\
& \leq \displaystyle{\sum_{\mid l-l'\mid\leq 4}}\parallel \parallel \left(S_{l'-1}v\right)^{\wedge}\parallel_{L^{1}(\mathbb{R}^{d})} \parallel \left(\phi_{l'}\hat{w}\right)\parallel_{L^{\vec{p}}(\mathbb{R}^{d})}\parallel_{L^{1}(I)}\\
&\leq\displaystyle{\sum_{\mid l-l'\mid\leq 4}}
\left(\displaystyle{\sum_{l''\leq l'}}2^{l''(s+1)))} \parallel \phi_{l''}\hat{v}\parallel_{L^{\infty}(I;L^{\vec{p}}(\mathbb{R}^{d}))}\right) \parallel \phi_{l'}\hat{w}\parallel_{L^{1}(I;L^{\vec{p}}(\mathbb{R}^{d}))}.
\end{split}    
\end{equation*}
Then,
\begin{equation*}
\parallel I^{l}_{1}\parallel_{L^{1}\left(I;L^{\vec{p}}(\mathbb{R}^{d})\right)}\leq 
\displaystyle{\sum_{\mid l-l'\mid\leq 4}}
\left( \displaystyle{\sum_{l''\leq l'}}2^{l''}2^{sl''}
\parallel \phi_{l''}\hat{v}\parallel_{L^{\infty}(I;L^{\vec{p}}(\mathbb{R}^{d}))}\right) \parallel \phi_{l'}\hat{w}\parallel_{L^{1}\left(I;L^{\vec{p}}(\mathbb{R}^{d})\right)}.
\end{equation*}
\end{proof}
By H\"older's inequality for series we get
\begin{equation*}
\parallel I^{l}_{1}\parallel_{L^{1}\left(I;L^{\vec{p}}(\mathbb{R}^{d})\right)}\leq \displaystyle{\sum_{\mid l-l'\mid\leq 4}}\left(\displaystyle{\sum_{l''\leq l'}}2^{l''q'} \right)^{\frac{1}{q'}}\cdot\parallel \phi_{l'}\hat{w}\parallel_{L^{1}\left(I;L^{\vec{p}}(\mathbb{R}^{d})\right)} \parallel v\parallel_{{\cal L}^{\infty}\left(I;\dot{{\cal FB}}^{s}_{\vec{p},q}\left(\mathbb{R}^{d}\right)\right)}.
\end{equation*}
If we multiply by $2^{l(s+1)}$ in both sides we have
\begin{equation*}
\begin{split}
2^{l(s+1)}\parallel I^{l}_{1}\parallel_{L^{1}\left(I;L^{\vec{p}}(\mathbb{R}^{d})\right)}&\leq
\left(
\displaystyle{\sum_{l'}}2^{(l-l')(s+1)}\chi_{\{k;\mid k\mid\leq 4\}}\left(l-l'\right)2^{l'(s+2)}\parallel \phi_{l'}\hat{w}\parallel_{L^{1}(I;L^{\vec{p}}(\mathbb{R}^{d}))}
\right)\cdot \parallel v\parallel_{Z}\\
&=\left(a_{k}\ast b_{l'}\right)_{l}\cdot \parallel v\parallel_{Z},
\end{split}    
\end{equation*}
where $a_{k}=2^{k(s+1)}\chi_{\{k;\mid k\mid\leq 4\}}\left(k\right)$ and $b_{l'}=2^{l'\left(s+2\right)}\parallel \phi_{l'}\hat{w}\parallel_{L^{1}(I;L^{\vec{p}}(\mathbb{R}^{d}))}$. By Young's inequality for series, we obtain
\begin{equation*}
\begin{split}
\parallel 2^{l(s+1)}\parallel I^{l}_{1}\parallel_{L^{1}\left(I;L^{\vec{p}}(\mathbb{R}^{d}) \right)}&\leq   \parallel a_{k}\parallel_{l^{1}\left(\mathbb{Z}\right)}\cdot 
\parallel b_{l'}\parallel_{l^{q}\left(\mathbb{Z}\right)}\cdot \parallel v\parallel_{Z}\\
&\leq \parallel w\parallel_{{\cal L}^{1}\left(I;\dot{{\cal FB}}^{s}_{\vec{p},q}\right)}\cdot \parallel v\parallel_{Z}\\
&\leq \parallel w\parallel_{Z}\cdot \parallel v\parallel_{Z}.
\end{split}    
\end{equation*}
Estimates for $I^{l}_{2}$: Similar computations as for $I^{l}_{1}$ provides
\begin{equation*}
\parallel 2^{l(s+1)}\parallel I^{l}_{2}\parallel_{L^{1}\left(I;L^{\vec{p}}(\mathbb{R}^{d}) \right)}\leq \parallel w\parallel_{Z} \parallel v\parallel_{Z}.     
\end{equation*}
Estimates for the term $I^{l}_{3}$: Since $\parallel \phi_{l}\hat{u}\parallel_{L^{1}(\mathbb{R}^{d})}\leq 2^{l(s+1)}\parallel \phi_{l}\hat{u}\parallel_{L^{\vec{p}}(\mathbb{R}^{d})}$, we get
\begin{equation*}
\begin{split}
 \parallel I^{l}_{3}\parallel_{L^{1}\left(I;L^{\vec{p}}(\mathbb{R}^{d})\right)}&\leq 
 \displaystyle{\sum_{l'\geq l-2}}\left(
 \displaystyle{\sum_{\mid l''-l'\mid\leq 1}}2^{l''(s+1))}\parallel \phi_{l''}\hat{w}\parallel_{L^{\infty}\left(I;L^{\vec{p}}(\mathbb{R}^{d})\right)}
 \right) 
 \parallel \phi_{l'}\hat{v}\parallel_{L^{1}\left(I;L^{\vec{p}}(\mathbb{R}^{d})\right)}\\
 &\leq  \displaystyle{\sum_{l'\geq l-2}}\left(
 \displaystyle{\sum_{\mid l''-l'\mid\leq 1}}2^{l''}\cdot 2^{s l''}\parallel \phi_{l''}\hat{w}\parallel_{L^{\infty}\left(I;L^{\vec{p}}(\mathbb{R}^{d})\right)}
 \right) 
 \parallel \phi_{l'}\hat{v}\parallel_{L^{1}\left(I;L^{\vec{p}}(\mathbb{R}^{d})\right)}.
\end{split}    
\end{equation*}
By the H\"older's inequality for series, we have
\begin{equation*}
\begin{split}
  \parallel I^{l}_{3}\parallel_{L^{1}\left(I;L^{\vec{p}}(\mathbb{R}^{d})\right)}&\leq
  \displaystyle{\sum_{l'\geq l-2}}
  \left(\sum_{\mid l''-l'\mid\leq 1}2^{l''q'}\right)^{\frac{1}{q'}} \parallel \phi_{l'}\hat{v}\parallel_{L^{1}\left(I;L^{\vec{p}}(\mathbb{R}^{d})\right)}\cdot\parallel w\parallel_{{\cal L}^{\infty}\left(I;\dot{{\cal FB}}^{s}_{\vec{p},q}(\mathbb{R}^{d})\right)}\\
  &\leq   \displaystyle{\sum_{l'\geq l-2}} 2^{-l'(s+1)}\left( 2^{l'(s+2)}\parallel \phi_{l'}\hat{v}\parallel_{L^{1}\left(I;L^{\vec{p}}(\mathbb{R}^{d})\right)}\right) \parallel w\parallel_{Z}
\end{split}    
\end{equation*}
By multiplying both sides by $2^{l(s+1)}$, we obtain
\begin{equation*}
\begin{split}
 2^{l(s+1)}  \parallel I^{l}_{3}\parallel_{L^{1}\left(I;L^{\vec{p}}(\mathbb{R}^{d})\right)}&\leq \displaystyle{\sum_{l'}}
 2^{\left(l-l'\right)\left(s+1\right)}\chi_{\{k;k\leq 2\}]}(l-l')\cdot \left( 2^{l'(s+2)}\parallel \phi_{l'}\hat{v}\parallel_{L^{1}\left(I;L^{\vec{p}}(\mathbb{R}^{d})\right)}\right)\cdot \parallel w\parallel_{Z}\\
 &\leq \left(a_{k}\ast b_{l'}\right)_{l}\cdot \parallel w\parallel_{Z},
\end{split}    
\end{equation*}
where $a_{k}=2^{k(s+1)}\chi_{\{k;k\leq 2\}}(k)$ and $b_{l'}=2^{l'(s+2)}\parallel \phi_{l'}\hat{v}\parallel_{L^{1}(I;L^{\vec{p}}(\mathbb{R}^{d}))}$.
Applying the Young's inequality for series, we obtain
\begin{equation*}
\begin{split}
\parallel 2^{l(s+1)}  \parallel I^{l}_{3}\parallel_{L^{1}\left(I;L^{\vec{p}}(\mathbb{R}^{d})\right)} \parallel_{l^{q}\left(\mathbb{Z}\right)}&\leq 
\parallel a_{k}\parallel_{l^{1}(\mathbb{Z})}\cdot\parallel b_{l'}\parallel_{l^{q}(\mathbb{Z})}\cdot \parallel w\parallel_{Z}\\
&\leq \parallel v\parallel_{Z}\cdot \parallel w\parallel_{Z},
\end{split}    
\end{equation*}
where the sequence $(a_{k})_{k\in\mathbb{Z}}\in l^{1}(\mathbb{Z})$ if $s+1>0$. These computations are enough to conclude the proof.

\begin{section}{Proof of the Main Theorem}
\label{lastsection}

In this section we conclude the proof of the main theorems. We set the computations for the mixed-norm Besov-Lebesgue space, since the procedure for the mixed-norm Fourier-Besov-Lebesgue spaces are similar. Let us have in mind expressions (\ref{mild-formulation}), (\ref{mild-formulation-abstract}) and Lemma \ref{fixed-point-scheme}. We first observe that, by Lemma \ref{estimates-lin-besov} we have
\begin{equation*}
 \parallel z_{0}\parallel_{Z}\leq c\max{(1,\frac{1}{\nu})}\parallel u_{0}\parallel_{\dot{{\cal B}}^{\sigma}_{\vec{p},q}(\mathbb{R}^{d})},   
\end{equation*}
where $c> 0$ does not depend on $\nu$. By Lemma \ref{estimates-bilinear-besov} we have
\begin{equation*}
 \parallel B\parallel_{{\cal B}(Z)}\leq K_{0}\max{(1,\frac{1}{\nu})} =:K.   
\end{equation*}
First we take 
\begin{equation*}
0<\varepsilon<\frac{1}{4K}=\frac{1}{4K_{0}\max{(1,\frac{1}{\nu})}}=\frac{\min{(1,\nu)}}{4K_{0}}.    
\end{equation*}
Therefore, in order to get the global well-posedness result it is enough to consider $u_{0}\in \dot{{\cal B}}^{\sigma}_{\vec{p},q}(\mathbb{R}^{d})$ such that
\begin{equation*}
 \parallel u_{0}\parallel_{\dot{{\cal B}}^{\sigma}_{\vec{p},q}(\mathbb{R}^{d})}\leq\frac{\varepsilon}{c\max{(1,\frac{1}{\nu})}}=\frac{\varepsilon}{c}\min{(1,\nu)}.   
\end{equation*}
This conclude the proof of item (i.) from the first main theorem. For the proof of the second theorem we proceed similarly. Finally, the procedure for the proof of the weak continuity from $[0,\infty)$ toward ${\cal S}'(\mathbb{R}^{d})$ is well-known and we omit it here.

\end{section}



\begin{thebibliography}{99}
	

\bibitem{Almeida-Ferreira-Lima-uniform-2017}{
Almeida, M.F.d., Ferreira, L.C.F. and Lima, L.S.M., Uniform global well-posedness of the Navier-Stokes-Coriolis system in a new critical space, Math. Z., 287, (2017), 735-750. 
https://doi.org/10.1007/s00209-017-1843-x
}


\bibitem{Benek-Panzone-mixed-Lp-1961}{
Benedek, A. and Panzone, R., The space $L^{p}$, with mixed norm, Duke Math. J., 28(1), (1961), 301-324.
}

\bibitem{Bony-symbolique-1981}{
Bony, J. M., Calcul symbolique et propagation des singularit\'es pour les \'equations aux d\'eriv\'ees partielles non lin\'eaires, In Annales scientifiques de l'\'Ecole normale sup\'erieure, 14(2), (1981), 209-246).
}

\bibitem{Bourgain-Pavlovic-2008}{
Bourgain, J. and Pavlovi\'c, N., Ill-posedness of the Navier-Stokes equations in a critical space in 3D, Journal of Functional Analysis, 255(9), (2008), 2233-2247.
}

\bibitem{Cannone-Wu-2012}{
Cannone M. and Wu G., Global well-posedness for Navier-Stokes equations in critical Fourier-Herz spaces, Nonlinear Analysis: Theory, Methods and Applications, 75(9), (2012) 3754-3760. https://doi.org/10.1016/j.na.2012.01.029.
}
\bibitem{Cannone-NS-generalization-1997}{
Cannone, M., A generalization of a theorem by Kato on Navier-Stokes equations, Revista Matem\'atica Iberoamericana, 13(3), (1997), 515-541.
}
\bibitem{Cleanthous-Georgiadis-Nielsen-mixed-norm-2016}{
Cleanthous, G., Georgiadis, A. and Nielsen, M., Spaces of distributions with mixed Lebesgue norms, In 15th Panhellenic Conference of Mathematical Analysis, (2016), 29-38.
}

\bibitem{Cleanthous-Georgiadis-Nielsen-discrete-2017}{
Cleanthous, G., Georgiadis, A. and Nielsen, M.,  Discrete decomposition of homogeneous mixed-norm Besov spaces, In Functional Analysis, Harmonic Analysis, and Image Processing: A collection of Papers in Honor of Bj\"orn Jawerth, American Mathematical Society, (2017), 167-184.
}

\bibitem{Bahouri-Chemin-Danchin-Fourier-PDE-2011}{
Bahouri, H., Chemin, J.-Y. and Danchin, R., Fourier Analysis and Nonlinear Partial Differential Equations, Series Title: Grundlehren der mathematischen Wissenschaften, Springer Berlin, Heidelberg, (2011). https://doi.org/10.1007/978-3-642-16830-7
}



\bibitem{Folland-real-1999}{
Folland, G. B., Real analysis: modern techniques and their applications, John Wiley and Sons,  (Vol. 40), (1999).
}
\bibitem{Furioli-Lemarie-Terrano-limite-2000}{
Furioli, G., Lemari\'e-Rieusset, P. G., and Terraneo E., Unicit\'e dans $L^{3}(\mathbb{R}^{3})$ et d'autres espaces fonctionnels limites pour Navier-Stokes, Revista Matem\'atica Iberoamericana, 16(3), (2000), 605-667.
}


\bibitem{Fang-Han-Hieber-Global-Navier-Stokes-Coriolis-Fourier-Besov-2015}
{\small Fang, D., Han, B. and Hieber, M., Global Existence Results for the Navier-Stokes Equations in the Rotational Framework in Fourier-Besov Spaces, In: Arendt, W., Chill, R., Tomilov, Y. (eds) Operator Semigroups Meet Complex Analysis, Harmonic Analysis and Mathematical Physics. Operator Theory: Advances and Applications, vol 250. Birkh\"auser, Cham. (2015).
}
 
\bibitem{Georgiadis-Nielsen-Pseudo-2016}{
Georgiadis, A. G. and Nielsen, M., Pseudodifferential operators on mixed-norm Besov and Triebel-Lizorkin spaces, Mathematische Nachrichten, 289(16), (2016), 2019-2036.
}



\bibitem{Iwabuchi-Takada-2014} {\small
Iwabuchi, T. and Takada, R., Global well-posedness and ill-posedness for the Navier-Stokes equations with the Coriolis force in function spaces of Besov type, Journal of Functional Analysis, 267.5 (2014), 1321-1337.
}



\bibitem{Khai-Tri-Sobolev-Lorentz-Navier-Stokes-2014}{
Khai, D.Q. and Tri, N.M., Solutions in mixed-norm Sobolev-Lorentz spaces to the initial value problem for the Navier-Stokes equations, Journal of Mathematical Analysis and Applications, 417(2), (2014), 819-833.
}

\bibitem{Kato-strong-1984}{
Kato, T., Strong $L^{p}$-solutions of the Navier-Stokes equation in $\mathbb{R}^{m}$, with applications to weak solutions, Mathematische Zeitschrift, 187, (1984), 471-480.
}

\bibitem{Konieczny-Yoneda-2011}{
Konieczny, P. and Yoneda, T., On dispersive effect of the Coriolis force for the stationary Navier-Stokes equations, Journal of Differential Equations, 250.10 (2011), 3859-3873.
}

\bibitem{Koch-Tataru-2001}{
Koch, H. and Tataru, D., Well-posedness for the Navier-Stokes equations, Advances in Mathematics, 157(1), (2001), 22-35.
}
\bibitem{Kozono-Yamazaki-semilinear-1994}{
Kozono, H. and Yamazaki, M., Semilinear heat equations and the Navier-Stokes equation with distributions in new function spaces as initial data, Communications in Partial Differential Equations, 19(5-6), (1994), 959-1014.
}

\bibitem{Lein-Lin-2011}{
Lei, Z. and Lin, F. H., Global mild solutions of Navier-Stokes equations, Communications on Pure and Applied Mathematics, 64(9), (2011), 1297-1304.
}
\bibitem{Lemarie-NS-century-2018}{
Lemari\'e-Rieusset, P. G., The Navier-Stokes problem in the 21st century, (2018), Chapman and Hall/CRC.
}




\bibitem{Paicu-anisotrope-2005}{
Paicu, M., \'Equation anisotrope de Navier-Stokes dans des espaces critiques, Revista Matem\'atica Iberoamericana, 21(1), (2005), 179-235.
}

\bibitem{Paicu-Zhang-global-2011}{
Paicu, M. and Zhang, P., Global Solutions to the 3-D Incompressible Anisotropic Navier-Stokes System in the Critical Spaces, Commun. Math. Phys., 307, (2011), 713-759.
https://doi.org/10.1007/s00220-011-1350-6
}

\bibitem{Planchon-NS-1996}{
Planchon, F., Global strong solutions in Sobolev or Lebesgue spaces to the incompressible Navier-Stokes equations in $\mathbb {R}^3$, Annales de l'I.H.P. Analyse non lin\'eaire, 13(3), (1996), 319-336.
}
\bibitem{Phan-mixed-2020}{
Phan, T., Well-posedness for the Navier-Stokes equations in critical mixed-norm Lebesgue spaces, Journal of Evolution Equations, 20(2), (2020), 553-576.
}

\bibitem{Planchon-NS-1998}{
Planchon, F., Asymptotic behavior of global solutions to the Navier-Stokes equations in $\mathbb R^ 3$, Revista Matem\'atica Iberoamericana, 14(1), (1998), 71-93.
}




\bibitem{Sawano-Besov-2018}{
Sawano, Y., Theory of Besov spaces, Singapore: Springer, (Vol. 56), (2018).
}
\bibitem{Sawano-homogeneous-2020}{
Sawano, Y., Homogeneous Besov spaces, Kyoto J. Math., 60(1), (2020), 1-43. https://doi.org/10.1215/21562261-2019-0038 
}










\end{thebibliography}
\end{document}